\documentclass[12pt,a4paper]{amsart}
\usepackage[T1]{fontenc}
\usepackage{lmodern}
\usepackage[english]{babel}
\usepackage[margin=1.2in]{geometry}
\usepackage{amsmath,amssymb,amsthm,mathtools,stmaryrd,mathrsfs}
\usepackage{microtype}
\usepackage{indentfirst}
\usepackage{needspace}
\usepackage{enumitem}
\usepackage[numbers,sort&compress]{natbib}
\usepackage[colorlinks=true,linkcolor=black,citecolor=black,urlcolor=black,
  bookmarksnumbered=true,hypertexnames=false]{hyperref}
\usepackage{bookmark}
\allowdisplaybreaks[2]
\numberwithin{equation}{section}
\newtheorem{theorem}{Theorem}[section]
\newtheorem{proposition}[theorem]{Proposition}
\newtheorem{lemma}[theorem]{Lemma}
\newtheorem{corollary}[theorem]{Corollary}
\theoremstyle{definition}
\newtheorem{definition}[theorem]{Definition}

\newtheorem{example}[theorem]{Example}
\theoremstyle{remark}
\newtheorem{remark}[theorem]{Remark}

\newcommand{\Zt}{\mathbb Z_2}
\newcommand{\R}{\mathbb R}

\newcommand{\Mass}{\mathbf M}
\newcommand{\Flat}{\mathcal F_r}
\newcommand{\dr}{\partial_r}
\newcommand{\II}{\mathbf I}
\newcommand{\IC}{\mathcal I}
\newcommand{\ZC}{\mathcal Z}
\newcommand{\VC}{\mathcal V}

\newcommand{\AMset}{\mathscr A}

\newcommand{\MC}{\llbracket M\rrbracket}
\newcommand{\res}{\mathbin{\llcorner}}
\newcommand{\HH}{\mathcal H}
\newcommand{\spt}{\operatorname{spt}}
\newcommand{\Int}{\operatorname{Int}}
\newcommand{\Vol}{\operatorname{Vol}}

\newcommand{\dist}{\operatorname{dist}}
\newcommand{\Lip}{\operatorname{Lip}}
\newcommand{\diver}{\operatorname{div}}
\newcommand{\SFF}{\mathrm{II}}
\newcommand{\dV}{\mathbf d_V}
\newcommand{\eps}{\varepsilon}
\newcommand{\An}{\operatorname{An}}
\newcommand{\id}{\operatorname{id}}
\hypersetup{
  pdftitle={Free-Boundary Minimal Varifolds I: Zero Boundary Mass},
  pdfauthor={},
  pdfsubject={Free-boundary stationary varifolds, almost minimizing property, boundary trace, zero boundary mass}
}
\title{Free-Boundary Minimal Varifolds I: Zero Boundary Mass}
\author{Tengyao Zhang}
\address{Yau Mathematical Sciences Center, Tsinghua University, Beijing, 100084, P.R. China}
\email{zhangtengyao25@mails.tsinghua.edu.cn}

\author{Zhifei Zhu}
\address{Yau Mathematical Sciences Center, Tsinghua University, Beijing, 100084, P.R. China}
\email{zhifeizhu@tsinghua.edu.cn}

\date{}
\begin{document}
\begin{abstract}
Let $V$ be a $k$-varifold in a compact Riemannian manifold $M^n$ with smooth nonempty boundary, where $2\le k<n$.
Assume that $V$ is stationary with free boundary and almost minimizing modulo two in every sufficiently small relative annulus about each point of $M$.
We prove that $V$ is integral and that the trace of the boundary second fundamental form on $T_xV$ vanishes at $\|V\|$-almost every boundary point.
On all of $M$, the first variation is $\delta V=\nu\,\sigma$, where $\nu$ is the outward unit normal and $\sigma$ is a nonnegative measure on $\partial M$ singular with respect to $\|V\|$.
A strict sign for the trace on all tangent $k$-planes at each boundary point forces zero boundary mass.
For hypersurfaces in ambient dimensions three through seven, strict mean convexity gives smoothness and proper embeddedness.
In every bounded smooth Euclidean domain $\Omega$, we obtain $\|V\|(\Omega)\ge c_{\Omega,k}\|V\|(\overline\Omega)$ with $c_{\Omega,k}>0$.
For smooth connected free-boundary minimal submanifolds in $\overline\Omega$, assuming a volume bound and a common annular almost minimizing condition at a fixed scale, then their first Dirichlet eigenvalues uniformly bounded above.
\end{abstract}
\maketitle

\setcounter{tocdepth}{2}
\begingroup
\small
\tableofcontents
\endgroup
\bigskip

\newpage
\section{Introduction}
\label{sec:introduction}

Let $(M^n,g)$ be a smooth compact connected Riemannian manifold of
dimension $n\ge3$ with nonempty smooth boundary, and let
$2\le k\le n-1$. A $k$-varifold $V$ in $M$ is said to be
\emph{stationary with free boundary} if its first variation vanishes on
every smooth vector field tangent to $\partial M$. For a smooth
properly embedded submanifold, this condition is equivalent to
minimality together with orthogonal contact along the boundary, and it
is the condition satisfied by the varifolds that arise from min--max
constructions in the space of relative cycles. For a general varifold
it is considerably weaker. Free-boundary stationarity allows both
nonzero normal first variation and positive boundary mass.
For instance, the $(n-1)$-varifold $|\partial M|$ is
stationary with free boundary for every metric $g$. Varifolds with
mass on the boundary also arise as limits of smooth free-boundary
minimal surfaces: in the Euclidean unit three-ball there are sequences
of smooth embedded free-boundary minimal surfaces that converge as
varifolds to the boundary sphere \cite{KS2024,KZ2026}; see
Subsection~\ref{sec:reference-comparison}.

We study varifolds that are \emph{almost minimizing modulo two in small
relative annuli}. This means that for every $p\in M$ there is a radius
$r_p>0$ such that $V$ is almost minimizing modulo two in every annulus
\[
 \An(p;s,r)=M\cap\bigl(B_r(p)\setminus\overline B_s(p)\bigr),
 \qquad 0<s<r<r_p,
\]
where the balls are taken in a fixed smooth extension of $M$.
For each annulus and each $\varepsilon>0$, Definition~\ref{def:almostminimizing}
requires a relative mod-two cycle approximating $V$ and a $\delta>0$
with the following property. Any deformation supported in the annulus
that lowers the mass of this cycle by more than $\varepsilon$ must
take a relative flat step larger than $\delta$ or raise the mass
above its initial value by more than $\delta$.

We prove that every
such $k$-varifold that is stationary with free boundary is integral, and that at
$\|V\|$-almost every point $x\in\partial M$ the trace of the second
fundamental form of $\partial M$ over the approximate tangent plane
$T_xV$ vanishes. For hypersurfaces this trace
is the mean curvature $H_{\partial M}$. In higher codimension it
depends on the tangent $k$-plane at the point of contact. Consequently,
a strict sign of the trace on all $k$-planes tangent to $\partial M$
excludes boundary mass. We also
show that $\delta V=\nu\,\sigma$ on all of $M$, where $\nu$ is the
outward unit normal and $\sigma$ is a nonnegative measure on $\partial M$
singular with respect to $\|V\|$.

This paper is the first in a series of two. It is concerned with the
structure of a given varifold: integrality, the first variation on
arbitrary vector fields, and the boundary trace identity. It contains
no existence statement. The second paper \cite{PartII}, referred to
below as Part~II, constructs varifolds to which the present results
apply. For every $2\le k\le n-1$, and under the assumption that
$H_j(M,\partial M;\Zt)=0$ for $1\le j<k$, it produces a nonzero
$k$-varifold that is stationary with free boundary and almost
minimizing modulo two in small relative annuli, and it bounds the mass
of this varifold in terms of the relative filling functions of $M$
\cite{PartII}. This is a free-boundary analogue of the
estimate of Nabutovsky and Rotman for the least mass of a nonzero
stationary integral $k$-varifold on a closed manifold \cite{NR2006}.
In the one-dimensional free-boundary problem, Beach, Contreras-Peruyero,
Griffin, Rotman, and Searle give an upper bound for the length of a
shortest orthogonal geodesic chord in a Riemannian three-ball with convex
boundary \cite{BCGRS2025}.

Throughout the paper, $(M^n,g)$ is as above. Chains and homology groups have coefficients in $\Zt$ unless otherwise stated. Stationarity with free boundary and almost minimization in small relative annuli are defined in Section~\ref{sec:notation}.

\subsection{Main results}

Let $\nu$ be the outward unit normal to $\partial M$, and let $\nabla$ be the Levi-Civita connection of $g$. We use the convention
\begin{equation}\label{eq:boundaryform}
 \SFF_{\partial M}(v,w)=g(\nabla_v\nu,w),
 \quad v,w\in T_x\partial M,
\end{equation}
for which the boundary of the Euclidean unit ball has positive second fundamental form. The scalar mean curvature of the boundary is
\begin{equation}\label{eq:boundarymean}
 H_{\partial M}=\operatorname{tr}_{T\partial M}\SFF_{\partial M},
\end{equation}
and we call the boundary strictly mean convex when $H_{\partial M}>0$ everywhere.

\begin{definition}
\label{def:kflat}
An integral $k$-varifold $V$ in $M$ is called $k$-flat if
\begin{equation}\label{eq:kflat}
 \operatorname{tr}_{T_xV}\SFF_{\partial M}=0
 \quad\text{for }\|V\|\text{-almost every }x\in\partial M.
\end{equation}
Here the approximate tangent plane $T_xV$ lies in $T_x\partial M$ at
$\|V\|$-almost every boundary point, by rectifiability, and the trace is
$\sum_{a=1}^k\SFF_{\partial M}(e_a,e_a)$ for an orthonormal basis
$e_1,\ldots,e_k$ of $T_xV$.
\end{definition}

\begin{definition}
\label{def:zeroboundarymass}
A $k$-varifold $V$ in $M$ has zero boundary mass if $\|V\|(\partial M)=0$.
\end{definition}

Every integral $k$-varifold with zero boundary mass is $k$-flat. For a smooth embedded submanifold $\Sigma\subset M$, proper embeddedness means $\Sigma\cap\partial M=\partial\Sigma$; it implies that $|\Sigma|$ has zero boundary mass. The converse need not hold, since zero boundary mass allows interior contact with $\partial M$ on a set of zero $k$-dimensional area.

Our main result is the following.

\begin{theorem}
\label{thm:kflat}
Let $2\le k\le n-1$. Suppose that a $k$-varifold $V$ in $M$ is stationary
with free boundary and almost minimizing modulo two in small relative annuli.
Then $V$ is integral and $k$-flat. There is a unique nonnegative Radon
measure $\sigma$ on $\partial M$, singular with respect to $\|V\|$, such that
\[
 \delta V(X)=\int_{\partial M}g(X,\nu)\,d\sigma
 \qquad\text{for every smooth vector field $X$ on $M$.}
\]
\end{theorem}

To prove Theorem~\ref{thm:kflat}, we first assume almost minimization
in a fixed relatively open set $U$. In this setting, almost minimization
itself implies free-boundary stationarity in $U$.

\begin{theorem}
\label{thm:local}
Let $2\le k\le n-1$, let $U\subset M$ be relatively open, and let $V$ be
a $k$-varifold in $M$ that is almost minimizing modulo two in $U$.
Then the following hold.
\begin{enumerate}[label=\textup{(\roman*)},leftmargin=2.2em]
\item The varifold $V$ is integral in $U$.
\item There is a unique nonnegative Radon measure $\sigma$ on
$U\cap\partial M$ such that
\[
 \delta V(X)=\int_{U\cap\partial M}g(X,\nu)\,d\sigma
\]
for every smooth vector field $X$ compactly supported in $U$.
\item The identity $\operatorname{tr}_{T_xV}\SFF_{\partial M}=0$ holds
for $\|V\|$-almost every $x\in U\cap\partial M$, and the measure
$\sigma$ in \textup{(ii)} is singular with respect to $\|V\|$ on $U$.
\end{enumerate}
\end{theorem}

Parts \textup{(i)} and \textup{(ii)} are
Proposition~\ref{loc:prop:integral}, and part \textup{(iii)} is
Proposition~\ref{kp:local-trace}. For a smooth properly embedded
free-boundary minimal submanifold $\Sigma$, the measure $\sigma$ is the
$(k-1)$-dimensional measure of $\partial\Sigma$, by
\eqref{eq:smoothvariation}. In general, the density of $\sigma$ with
respect to $\|V\|$ on $\partial M$ equals the trace in \eqref{eq:kflat}
(Theorem~\ref{kp:thm:nonnegative-trace}), so the two assertions of
\textup{(iii)} are equivalent. In Subsection~\ref{sec:consequences}
we apply Theorem~\ref{thm:local} separately to the annuli covering each
punctured ball. The measures from these annuli agree on overlaps.
Global free-boundary stationarity gives the monotonicity estimate at the
centers of the cover. With this estimate, a cutoff argument extends the
first-variation formula across the centers, yielding Theorem~\ref{thm:kflat}.

The trace identity has the following consequences for the boundary
mass.

\begin{corollary}
\label{cor:zeroboundary}
Assume that, at every $x\in\partial M$, one of the following alternatives holds:
\begin{equation}\label{eq:stricttrace}
 \begin{aligned}
 &\operatorname{tr}_P\SFF_{\partial M}>0
 &&\text{for every $k$-plane }P\subset T_x\partial M,\\
 \text{or}\quad
 &\operatorname{tr}_P\SFF_{\partial M}<0
 &&\text{for every $k$-plane }P\subset T_x\partial M.
 \end{aligned}
\end{equation}
Then every varifold $V$ as in Theorem~\ref{thm:kflat} has zero boundary mass. The sign in \eqref{eq:stricttrace} is constant on each connected component of $\partial M$ and may differ between components.
\end{corollary}

For hypersurfaces in ambient dimensions three through seven, regularity
theory gives a stronger conclusion under strict mean convexity.

\begin{corollary}
\label{cor:hypersurface}
Suppose $3\le n\le7$, $k=n-1$, and $H_{\partial M}>0$ on $\partial M$. Every nonzero $(n-1)$-varifold that is stationary with free boundary and almost minimizing modulo two in small relative annuli has a representation
\begin{equation}\label{eq:hypersurfacerepresentation}
 V=\sum_{i=1}^{N}a_i|\Sigma_i|,
\end{equation}
where $N$ and the multiplicities $a_i$ are positive integers, and each $\Sigma_i$ is a smooth, compact, connected, properly embedded free-boundary minimal hypersurface, possibly with empty boundary. In particular, $\|V\|(\partial M)=0$.
\end{corollary}

In a Euclidean domain, the trace identity gives a lower bound for the
fraction of mass in the interior, without any boundary curvature assumption.

\begin{corollary}
\label{cor:euclidean-interior}
Let $\Omega\subset\R^n$ be a bounded connected domain with smooth
boundary, and equip $M=\overline\Omega$ with the Euclidean metric.
For each $2\le k\le n-1$, there is a constant $c_{\Omega,k}>0$ such
that every $k$-varifold $V$ in $M$ that is stationary with free
boundary and almost minimizing modulo two in small relative annuli satisfies
\[
 \|V\|(\Omega)\ge c_{\Omega,k}\,\|V\|(\overline\Omega).
\]
\end{corollary}

For smooth free-boundary minimal submanifolds of bounded area, almost
minimization at a fixed scale also gives a lower bound for the depth.
Let $\Omega\subset\R^n$ be a bounded connected domain with smooth
boundary, let $2\le k\le n-1$, and equip $M=\overline\Omega$ with the
Euclidean metric. In Definition~\ref{def:almostminimizing} we take
$\R^n$ with the Euclidean metric as the extension, so that the relative
annuli are $\An(p;s,r)=\{x\in\overline\Omega:s<|x-p|<r\}$. Fix an
integer $N\ge1$ and a radius $R_*>0$. A $k$-varifold in
$\overline\Omega$ has \emph{property $\mathcal P_N$ at scale $R_*$} if,
for every $p\in\overline\Omega$ and all radii satisfying
\begin{equation}\label{eq:PN}
 0<s_j<r_j<R_*\quad(1\le j\le N),
 \quad
 s_j>12\,r_{j+1}\quad(1\le j<N),
\end{equation}
it is almost minimizing modulo two in at least one of the relative
annuli $\An(p;s_j,r_j)$, $1\le j\le N$. Almost minimization arises in
this form in min--max constructions, through the combinatorial argument
of Pitts \cite{Pitts1981}. In contrast to almost minimization in small
relative annuli, where the radius $r_p$ depends on the varifold,
property $\mathcal P_N$ at a fixed scale passes to varifold limits;
compare \cite[Lemma~3.9]{PartII}. The factor $12$ follows the convention
of \cite[Section~3]{PartII} and plays no role in the proofs.

By a free-boundary minimal submanifold of $\overline\Omega$ we mean a
smooth compact embedded $k$-dimensional submanifold
$\Sigma\subset\overline\Omega$ with
$\partial\Sigma\subset\partial\Omega$ that has zero mean curvature and
meets $\partial\Omega$ orthogonally along $\partial\Sigma$. Proper
embeddedness is not assumed, so the interior of $\Sigma$ may touch
$\partial\Omega$. The depth of $\Sigma$ is
\[
 d_\Sigma=\max_{x\in\Sigma}\dist(x,\partial\Omega).
\]
When $\Sigma$ is connected and $\partial\Sigma\ne\emptyset$,
$\lambda_1^D(\Sigma)$ denotes the first eigenvalue of the Laplacian of
the induced metric with Dirichlet condition on $\partial\Sigma$; no
condition is imposed where the interior of $\Sigma$ touches
$\partial\Omega$.

\begin{corollary}[Depth and Dirichlet eigenvalue]
\label{cor:dirichlet}
Let $\Omega\subset\R^n$ be a bounded connected domain with smooth
boundary, and let $2\le k\le n-1$. For every $\Lambda>0$, every integer
$N\ge1$, and every $R_*>0$ there is a constant
$d_*=d_*(\Omega,k,\Lambda,N,R_*)>0$ with the following property. Let
$\Sigma\subset\overline\Omega$ be a nonempty connected $k$-dimensional
free-boundary minimal submanifold such that $\HH^k(\Sigma)\le\Lambda$
and the multiplicity-one varifold $|\Sigma|$ has property
$\mathcal P_N$ at scale $R_*$. Then the boundary $\partial\Sigma$ is
nonempty,
\begin{equation}\label{eq:depth}
 d_\Sigma\ge d_*,
\end{equation}
and
\begin{equation}\label{eq:dirichlet}
 \lambda_1^D(\Sigma)
 \le\frac{(k+1)(k+2)}{2\omega_k}\,
     \frac{\HH^k(\Sigma)}{d_\Sigma^{\,k+2}}
 \le\frac{(k+1)(k+2)}{2\omega_k}\,
     \frac{\Lambda}{d_*^{\,k+2}},
\end{equation}
where $\omega_k$ is the volume of the Euclidean unit $k$-ball. No
assumption is made on the curvature of $\partial\Omega$.
\end{corollary}

The depth bound \eqref{eq:depth} is where the trace identity enters. If
it failed, a sequence of such submanifolds would converge to a nonzero
varifold that is stationary with free boundary, almost minimizing in
small relative annuli, and supported in $\partial\Omega$, contrary to
Corollary~\ref{cor:euclidean-interior}. The eigenvalue bound then
follows from interior monotonicity and an explicit test function; its
first inequality involves only the depth of $\Sigma$ itself. An area bound alone gives neither a uniform positive lower bound for the depth
nor a uniform upper bound for the first Dirichlet eigenvalue. For the free-boundary minimal
surfaces of \cite{KS2024} in the unit three-ball, the areas are
bounded, the depth tends to zero, and the first Dirichlet eigenvalue
tends to infinity; see the end of Subsection~\ref{sec:consequences}.

For Corollary~\ref{cor:euclidean-interior}, the trace identity and
Euclidean first variation give the interior mass bound in every
codimension \cite{Allard1972}. Harvey and Lawson give a related
convex-function argument for currents
\cite[Theorem~4.10 and Corollary~4.11]{HL2013}. For hypersurfaces,
the trace identity gives zero boundary mass even when $H_{\partial M}<0$.
When $H_{\partial M}>0$ and $3\le n\le7$, regularity
\cite[Theorem~4.6]{GLWZ2021} and the maximum principle
(Proposition~\ref{prop:barrier}) also give properness. For properness
under a generic metric, see \cite{SWZ2024}, and for free-boundary
min--max existence in higher dimensions, see
\cite[Theorem~1.8]{BLPS2026}.

The examples in Subsection~\ref{sec:examples} explain the role of the
hypotheses. The boundary sphere in Example~\ref{ex:spherical-boundary}
is locally minimizing among absolute cycles and has negative boundary
trace. Wang's obstacle problem likewise permits an outward mean
curvature vector on the contact set \cite{Wang2022}. In
Proposition~\ref{prop:boundary-concentration}, a varifold
almost minimizing among relative cycles carries all its mass on a
totally geodesic boundary.

\subsection{Related work}
\label{sec:reference-comparison}

In the Euclidean unit three-ball, genus-zero free-boundary minimal
surfaces maximizing the first normalized Steklov eigenvalue converge
as varifolds to the boundary sphere as the number of boundary components
tends to infinity \cite[Corollary~1.5]{KS2024}. Kapouleas and Zou obtained the same
limit by gluing \cite[Theorem~1.1]{KZ2026}. Boundary concentration also
occurs for Neumann Allen--Cahn solutions in the unit ball \cite{MNW2007}
and in strictly mean-convex Euclidean domains \cite{MW2007}, and for
Ginzburg--Landau solutions in codimension two \cite{Cheng2020}.

In relative min--max theory, canonical representatives have zero boundary
mass, but their varifold limits can charge $\partial M$. In dimensions
$3\le n\le7$, Li and Zhou proved smoothness even at interior points of
contact with $\partial M$ \cite{LZ2021}. In a strictly mean-convex
three-ball, the Simon--Smith construction gives an embedded free-boundary
minimal disk \cite[Theorem~1.2]{SSWZ2026}.

For Neumann Allen--Cahn limits, Li, Parise, and Sarnataro proved
integrality up to the boundary \cite[Theorem~3.3]{LPS2024}. For stable
solutions, the first variation is singular with respect to the limit
weight, and boundary mass is concentrated where $H_{\partial M}=0$
\cite[Theorem~4.8(3) and Remark~4.11]{LPS2026}.
Theorems~\ref{thm:kflat} and~\ref{thm:local} give the corresponding
results for relative almost minimizing varifolds in arbitrary codimension.

For stationary free-boundary varifolds, known results include regularity
near a multiplicity-one half-plane \cite{GJ1986} and boundary
monotonicity in arbitrary codimension \cite{GLZ2020}. Li and Zhou's
maximum principle excludes one-sided contact with a strongly $k$-convex
barrier meeting $\partial M$ orthogonally \cite{LZBarrier2021}. The trace
identity \eqref{eq:kflat} concerns the curvature of $\partial M$ itself
and eliminates the boundary curvature term in De~Masi's Euclidean
first-variation formula \cite[Corollary~4.7]{DeMasi2021}.

In capillary theory, the contact angle determines the weight of the
wetted area in the energy \cite{LZZ2025}, and separate first-variation
bounds hold for the interior and boundary portions \cite{DEGL2025}.
Here chains supported on $\partial M$ represent zero as relative cycles.

\subsection{Proof strategy and organization}
\label{sec:proof-overview}

Section~\ref{sec:notation} introduces relative chains and varifolds.
In Section~\ref{sec:local} we use relative almost minimization to
recover the normal component of the first variation. We replace the approximating
cycles by comparison currents whose limiting mass measure agrees with
that of the prescribed $V$, including any boundary mass. On each
comparison neighborhood the currents $R_i$ have weighted first variation
$\delta_{w_i}|R_i|=\nu\,w_i\sigma_i$, where $w_i$ is the weight and
$\sigma_i$ is a nonnegative measure on $\partial M$. As $w_i\to1$,
we obtain integrality of $V$ and the local identity $\delta V=\nu\,\sigma$
for a nonnegative boundary measure $\sigma$.

In Section~\ref{sec:kflat}, the locality of mean curvature identifies
the density of $\sigma$ with respect to $\|V\|$ with the boundary trace,
which is therefore nonnegative. To exclude positive trace, we combine an
inward deformation fixing the boundary pointwise with a small expansion
preserving the boundary. The expansion makes the inward area gain
dominate the cost of the thin transition collar, and quantitative
height and tilt estimates control the remaining errors. Minimality of
an auxiliary functional in the relative comparison class then yields a
contradiction. The section concludes with the geometric consequences
and examples showing why the hypotheses are needed.
Part~II applies these results to its min--max varifolds in
\cite[Corollaries~1.6 and~1.7]{PartII}.

	\subsection{Acknowledgments}
	The authors are supported in part by the NSFC Grant 12501065. Z.Z. is also supported in part by National Key R\&D Program of China 2024YFA1015300.
	
  \subsection{Statements and Declarations}
	Data sharing is not applicable to this article as no datasets were generated or analyzed during the current study. ChatGPT were used to assist with bibliographic searches, language editing, and initial drafts of portions of some proofs.  Every AI-assisted mathematical argument was independently reconstructed and checked by the authors, who reviewed and revised all retained material and take full responsibility for the manuscript. 

\section{Relative chains and varifolds}
\label{sec:notation}

We use integral flat chains with coefficients in $\Zt$ on $(M,\partial M)$ \cite{Federer1969}, equivalently the rectifiable relative-cycle model of Li--Zhou \cite[Section 3]{LZ2021} with mod-two coefficients \cite[Section 3]{GLWZ2021}.
Informally, a relative cycle is a chain whose boundary lies in $\partial M$, and chains supported in $\partial M$ count as zero.

\begin{definition}
\label{def:relativechains}
Let $\II_j(M;\Zt)$ denote integral $j$-chains supported in $M$, and put
\[
 \IC_j(M,\partial M;\Zt)
 =\II_j(M;\Zt)/\II_j(\partial M;\Zt),
 \quad \dr[T]=[\partial T].
\]
The kernel of $\dr$ is denoted by $\ZC_j(M,\partial M;\Zt)$. For a relative class $T$, define
\[
 \Mass(T)=\inf\{\Mass(\widetilde T):[\widetilde T]=T\}.
\]
For relative cycles, we also use the rectifiable model: finite-mass rectifiable $j$-chains whose boundaries are supported in $\partial M$, identified when their difference is supported in $\partial M$. Each such class has a unique canonical rectifiable representative with zero mass on $\partial M$. Its mass equals the displayed infimum over integral representatives. We omit brackets for relative classes when there is no ambiguity.
\end{definition}

The equivalence and approximation by integral representatives are given in \cite[Lemma 3.1]{GLWZ2021} and \cite[Proposition 2.4]{LMN2018}. A canonical representative may have infinite absolute boundary mass, so the infimum need not be attained. We interpret chain identities using integral representatives before passing to the quotient, and analytic limits using canonical representatives. For Borel $E\subset M$, $\|T\|(E)$ is the canonical representative's mass in $E$; its ordinary varifold is $|T|$.

The relative flat norm is
\begin{equation}\label{eq:flatnorm}
 \Flat(T)=\inf_{T=R+\dr Q}
 \bigl(\Mass(R)+\Mass(Q)\bigr),
\end{equation}
where $R$ and $Q$ are relative chains of dimensions $j$ and $j+1$. The flat and mass distances between two classes are $\Flat(T-S)$ and $\Mass(T-S)$, respectively. In particular,
\begin{equation}\label{eq:flatbasic}
 \Flat(T-S)\le\Mass(T-S),\quad
 \Flat(\dr Q)\le\Mass(Q).
\end{equation}
For a chain $T$ and a measurable set $E$, $T\res E$ denotes restriction. For a Lipschitz map $f$, $f_\#T$ denotes pushforward. The mod-two relative fundamental chain is $\MC$, and $[M,\partial M]\ne0$ denotes its homology class.

\begin{definition}
\label{def:stationarity}
Let $G_k(M)$ be the bundle of unoriented $k$-planes in $TM$. A $k$-varifold is a nonnegative Radon measure on $G_k(M)$; the space of these measures is denoted by $\VC_k(M)$. Its mass measure $\|V\|$ is the pushforward under the bundle projection. An integral varifold is represented by a countably $k$-rectifiable set with positive integer multiplicity almost everywhere.

Let $\mathfrak X_{\tan}(M)$ be the smooth vector fields on $M$ tangent to $\partial M$ along the boundary. The first variation of $V$ is
\begin{equation}\label{eq:firstvariation}
 \delta V(X)=\int_{G_k(M)}\diver_PX(x)\,dV(x,P),
 \quad
 \diver_PX=\sum_{a=1}^{k}g(\nabla_{e_a}X,e_a),
\end{equation}
where $e_1,\ldots,e_k$ is an orthonormal basis of $P$. The varifold is stationary with free boundary if $\delta V(X)=0$ for every $X\in\mathfrak X_{\tan}(M)$. A local statement in an open set $U\subset M$ uses fields compactly supported in $U$.
\end{definition}

For a smooth submanifold $\Sigma$ with outward conormal $\eta_\Sigma$, its first variation is
\begin{equation}\label{eq:smoothvariation}
 \delta|\Sigma|(X)
 =-\int_\Sigma g(\mathbf H_\Sigma,X)\,d\HH^k
  +\int_{\partial\Sigma}g(\eta_\Sigma,X)\,d\HH^{k-1},
\end{equation}
where $\mathbf H_\Sigma$ is the mean curvature vector and $\HH^j$ is $j$-dimensional Riemannian Hausdorff measure. A smooth free-boundary minimal submanifold has $\mathbf H_\Sigma=0$ and meets $\partial M$ orthogonally along its boundary.

We use the free-boundary version of Pitts' almost minimizing condition
\cite{Pitts1981} developed in \cite[Definitions~3.18 and~4.19]{LZ2021}.
For the mod-two hypersurface case, see \cite[Section~4]{GLWZ2021}.

\begin{definition}
\label{def:almostminimizing}
Fix a smooth Riemannian metric on the compact Grassmann bundle $G_k(M)$, and let $d_G$ be its distance. For $V,W\in\VC_k(M)$, define
\begin{equation}\label{eq:varifold-distance}
 \dV(V,W)=\sup_{\substack{\|f\|_\infty\le1\\\Lip_{d_G}(f)\le1}}
 \left|\int_{G_k(M)}f\,dV-\int_{G_k(M)}f\,dW\right|.
\end{equation}
The supremum is over real-valued functions on $G_k(M)$ whose supremum norm and Lipschitz constant with respect to $d_G$ are at most one. This distance induces weak convergence on sets of varifolds with bounded mass.
Given an open set $U\subset M$, $\eps,\delta>0$, and $\mathfrak d\in\{\Flat,\Mass\}$, let $\AMset(U;\eps,\delta;\mathfrak d)$ be the relative $k$-cycles $T$ such that every finite sequence $T=T_0,\ldots,T_N$ satisfying
\[
 \spt(T_i-T)\subset U,\quad
 \mathfrak d(T_i-T_{i-1})\le\delta,\quad
 \Mass(T_i)\le\Mass(T)+\delta
\]
also satisfies $\Mass(T_N)\ge\Mass(T)-\eps$. Supports are taken for canonical representatives.

A varifold $V$ is almost minimizing modulo two in $U$ if, for every $\eps>0$, there are $\delta>0$ and $T\in\AMset(U;\eps,\delta;\Flat)$ such that $\dV(|T|,V)<\eps$.

Choose a smooth extension of the metric across $\partial M$. For balls in this extension, define the relative annulus
\[
 \An(p;s,r)=M\cap\bigl(B_r(p)\setminus\overline B_s(p)\bigr),
 \quad0<s<r.
\]
The varifold is almost minimizing in small relative annuli if, for each $p\in M$, there is $r_p>0$ such that it is almost minimizing in every $\An(p;s,r)$ with $0<s<r<r_p$.
\end{definition}

In words, a cycle $T$ belongs to $\AMset(U;\eps,\delta;\mathfrak d)$ if its mass cannot be decreased by more than $\eps$ through a finite chain of modifications inside $U$, each of which is $\delta$-small and none of which raises the mass by more than $\delta$. A varifold is almost minimizing if it is approximated by such cycles with $\eps\to0$. The cycles and the steps are relative, whereas $|T|$ and $V$ are ordinary varifolds on $M$.

Canonical representatives have zero boundary mass, but their ordinary varifold limits may acquire it. Theorem~\ref{thm:kflat} addresses this issue.

We finally record the compactness and small-filling properties of relative cycles.

\begin{lemma}
\label{inp:compactness}
Let $1\le k\le n-1$. Relative cycles with a common mass bound are sequentially compact in the relative flat topology, and relative mass is lower semicontinuous. There are constants $\rho_M,C_M>0$ such that, if
$\Flat(T-S)<\rho_M$, then there is a relative $(k+1)$-chain $Q$ satisfying
\begin{equation}\label{inp:eq:small-filling}
 \dr Q=T-S,
 \quad\Mass(Q)\le C_M\Flat(T-S).
\end{equation}
\end{lemma}

Compactness and lower semicontinuity follow from
\cite[Section~2.3 and Proposition~2.4]{LMN2018}. In the compactness
argument, projection in a boundary collar gives equivalent integral
representatives with uniformly bounded boundary mass, to which
Federer--Fleming compactness applies.

For the small-filling estimate, apply Almgren's relative isoperimetric
choice \cite[Proposition~1.23]{Almgren1962} with mod-two coefficients,
as in \cite[Appendix~A, p.~957]{LMN2018}. For integral representatives
$\widetilde T,\widetilde S$ of the two cycles, it gives a mod-two chain
$\widetilde Q$ in $M$ of mass at most $C_M\Flat(T-S)$ such that
$\partial\widetilde Q-(\widetilde T-\widetilde S)$ is supported in
$\partial M$. Its relative class $Q=[\widetilde Q]$ therefore satisfies
\eqref{inp:eq:small-filling}.

\section{Comparison currents and integrality}
\label{sec:local}

Free-boundary stationarity alone does not imply rectifiability. In this section we use almost minimization to prove integrality and determine the first variation on arbitrary vector fields. The statement is local. It is the starting point for the trace identity of Section~\ref{sec:kflat}, and it gives the integrality of the min--max varifolds constructed in \cite{PartII}.

We use the outward normal $\nu$, extended along normal geodesics near $\partial M$, and the convention \eqref{eq:boundaryform}. The connection $\nabla$ and all unmarked norms refer to $g$. For a subspace $P\subset T_xM$, $\pi_P$ denotes orthogonal projection and $P^\perp$ its orthogonal complement. Balls are taken in the smooth ambient extension fixed in Section~\ref{sec:notation}, and intersected with $M$. A field is admissible if it is smooth, compactly supported in the open set under consideration, and tangent to $\partial M$. Its flow preserves $M$.

\begin{proposition}
\label{loc:prop:integral}
Let $2\le k\le n-1$, let $U\subset M$ be relatively open, and let $V$
be a $k$-varifold in $M$ that is almost minimizing modulo two in $U$.
Then $V$ is integral in $U$, and there is a unique nonnegative Radon
measure $\sigma$ on $U\cap\partial M$ such that
\[
 \delta V(X)=\int_{U\cap\partial M}g(X,\nu)\,d\sigma
\]
for every smooth vector field $X$ compactly supported in $U$.
\end{proposition}

By \cite[Proposition~3.21]{LZ2021}, almost minimization implies that
$V$ is stationary with free boundary in $U$. The same argument applies
with mod-two coefficients since a boundary-preserving flow acts on
relative mod-two cycles and is continuous in relative flat distance.

The proof has three steps. We replace the cycles approximating $V$ by comparison currents $R_i$ that minimize an auxiliary functional (Lemma~\ref{kp:lem:comparison-minimizers}). We then compute the first variation of these currents on arbitrary vector fields and estimate their mass in a thin collar of the boundary (Lemma~\ref{kp:lem:boundary-force}). Finally we pass to the limit (Lemma~\ref{kp:lem:limit-form}). The currents $R_i$ and the collar estimate are used again in Section~\ref{sec:kflat}.

All constants are local: $C$ is independent of sequence indices and
shrinking radii, while $C_i$ may depend on the fixed current indexed
by $i$. The notation $O(r^a)$ and $o(r^a)$ has its usual meaning
with these constants. All currents have coefficients in $\Zt$.

\subsection{Comparison currents}
\label{kp:subsec:comparison-currents}

Throughout this section, $V$ and $U$ are as in Proposition~\ref{loc:prop:integral}.
Almost minimization is formulated with relative classes and flat steps
(Definition~\ref{def:almostminimizing}). The comparison varifolds
$|R_i|$ constructed below have zero boundary mass, but their ordinary
varifold limit may carry mass on $\partial M$. We replace the cycles
approximating $V$ by minimizers of an auxiliary functional. This
functional controls finitely many test integrals at each stage and
thereby preserves the limiting mass measure.

Write $\|f\|_{C^m}=\max_{0\le q\le m}\sup_M|\nabla^qf|_g$,
where $\nabla^0f=f$. Embed $M$ smoothly in Euclidean space
\cite[Chapter 6]{Lee2013}, enumerate the coordinate monomials with
$\varphi_1=1$, and rescale so that $\|\varphi_j\|_{C^j}\le1$.
Their span remains uniformly dense in $C(M)$ by Stone--Weierstrass
\cite[Theorem 7.32]{Rudin1976}.

Since $\|\varphi_j\|_{C^m}\le1$ for $j\ge m$, separating the finitely many terms with $j<m$ gives constants $C_m\ge1$, independent of $i$ and with $C_0=C_1=1$, such that
\begin{equation}\label{kp:eq:test-bound}
 \frac{1}{i^2}\sum_{j=1}^i\|\varphi_j\|_{C^m}^2\le\frac{C_m}{i},
 \quad
 \frac{1}{i^2}\sum_{j=1}^i\|\varphi_j\|_{C^0}\|\varphi_j\|_{C^m}
 \le\frac{C_m}{i}.
\end{equation}
Set $\eps_i=i^{-4}$. Definition~\ref{def:almostminimizing} supplies relative cycles $\tau_i$, canonical representatives $T_i$, and $0<\delta_i\le1/i$ such that
\[
 \tau_i\in\AMset(U;\eps_i,\delta_i;\Flat),\quad |T_i|\longrightarrow V.
\]
Shrinking $\delta_i$ preserves this property.

Fix open sets $\Omega'\Subset\Omega\Subset U$, where $\Subset$ means compact containment in the relative topology of $M$, and put $K=\overline\Omega$.
Adapting the comparison class in \cite[proof of Proposition~5.3, Step~1]{LZ2021}, let $\mathcal C_i$ consist of the endpoints of finite sequences starting at $\tau_i$ and satisfying
\begin{equation}\label{kp:eq:reachable}
 \spt(\tau^{(j)}-\tau_i)\subset K,\quad
 \Flat(\tau^{(j+1)}-\tau^{(j)})\le\delta_i,\quad
 \Mass(\tau^{(j)})\le\Mass(T_i)+\delta_i.
\end{equation}
These conditions hold at every available index, and a sequence of length zero is allowed. Write $[S]$ for the relative class of a canonical current $S$. Almost minimizing gives
\begin{equation}\label{kp:eq:am-lower}
 \Mass(S)\ge\Mass(T_i)-\eps_i\quad([S]\in\mathcal C_i).
\end{equation}
For a finite nonnegative measure $\alpha$, define
\begin{equation}\label{kp:eq:auxiliary-functional}
 F_i(\alpha)=\alpha(M)+\frac{1}{i^2}\sum_{j=1}^i
 \left(\int\varphi_j\,d(\alpha-\|T_i\|)\right)^2,
 \quad F_i(S)=F_i(\|S\|).
\end{equation}

The role of the penalty terms is the following. A minimizer of mass in $\mathcal C_i$ would have a good first variation, but its mass measure could differ from $\|T_i\|$, and the limit of these measures could differ from $\|V\|$. A minimizer of $F_i$ instead has its first $i$ test integrals close to those of $T_i$; the price is that mass is replaced by a weighted mass, with a weight close to one.
Such quadratic terms occur in White's local minimax construction
\cite{White1994} and Marques--Neves's auxiliary mass functional
\cite[proof of Theorem 6.1 and proof of Proposition 6.2, Claim 2]{MN2021}.
Here the moving centers $\|T_i\|$ and increasing number of tests
preserve the prescribed limiting mass measure.

\begin{lemma}
\label{kp:lem:comparison-minimizers}
For all sufficiently large $i$, $F_i$ has a minimum among canonical currents with class in $\mathcal C_i$, attained at $R_i$, and
\begin{equation}\label{kp:eq:comparison-bounds}
 \begin{gathered}
 \Mass(T_i)-\eps_i\le\Mass(R_i)\le F_i(R_i)\le\Mass(T_i),\\
 \left|\int\varphi_j\,d(\|R_i\|-\|T_i\|)\right|\le\frac1i
 \quad(1\le j\le i).
 \end{gathered}
\end{equation}
Here $R_i$ is a canonical rectifiable current; finite mass of its absolute boundary is not asserted. Its associated $|R_i|$ is an ordinary integral varifold. In particular, $\|R_i\|\rightharpoonup\|V\|$. There are positive smooth functions $w_i\to1$ in every fixed $C^m$ norm such that $|R_i|$ is stationary with free boundary for $g_i=w_i^{2/k}g$ in $\Omega$.
\end{lemma}

\begin{proof}
We begin with the derivative of $F_i$ and a monotonicity property.
Choose $L$ with $\Mass(T_i)+\delta_i<L$ for every $i$. Differentiating $F_i$ at $\alpha$ in a signed measure direction $\beta$ gives $\int w_i[\alpha]\,d\beta$, where
\begin{equation}\label{kp:eq:weight}
 w_i[\alpha](x)=1+\frac{2}{i^2}\sum_{j=1}^i
 \left(\int\varphi_j\,d(\alpha-\|T_i\|)\right)\varphi_j(x).
\end{equation}
For $\alpha(M)\le L$, the integrals in \eqref{kp:eq:weight} are bounded by $2L\|\varphi_j\|_{C^0}$, so
\[
 \|w_i[\alpha]-1\|_{C^m}
 \le\frac{4L}{i^2}\sum_{j=1}^i
        \|\varphi_j\|_{C^0}\|\varphi_j\|_{C^m}
 \le\frac{4LC_m}{i}
\]
by \eqref{kp:eq:test-bound}.
Thus $1/2\le w_i[\alpha]\le3/2$ for large $i$, uniformly when $\alpha(M)\le L$. For $0\le\beta\le\alpha$ in this range, integration along the segment gives
\begin{equation}\label{kp:eq:measure-monotonicity}
 F_i(\alpha)-F_i(\beta)
 =\int_0^1\int w_i[\beta+t(\alpha-\beta)]\,d(\alpha-\beta)\,dt
 \ge\tfrac12(\alpha-\beta)(M)\ge0.
\end{equation}

Next we show that $F_i$ attains its minimum on $\mathcal C_i$.
Fix such an $i$ and take a canonical minimizing sequence with $[S_\ell]\in\mathcal C_i$. Lemma~\ref{inp:compactness} and weak compactness of measures give a canonical $R_i$ and a measure $\alpha$ such that
\[
 \Flat([S_\ell]-[R_i])\to0,\quad
 \|S_\ell\|\rightharpoonup\alpha,
 \quad \alpha(M)\le\Mass(T_i)+\delta_i.
\]
Lower semicontinuity gives $\|R_i\|\le\alpha$ in the interior and, since $R_i$ is canonical, on all of $M$. Relative flat convergence gives $R_i=T_i$ outside $K$ in the interior; canonicality then gives $\spt(R_i-T_i)\subset K$. For large $\ell$, appending $[R_i]$ to a sequence reaching $[S_\ell]$ satisfies \eqref{kp:eq:reachable}: its last flat step is less than $\delta_i$ and its mass is at most $\alpha(M)$. Hence $[R_i]\in\mathcal C_i$.

Passing to the canonical representative may lose mass, including mass
accumulated on $\partial M$. The quadratic penalties alone need not
decrease under such a loss; by \eqref{kp:eq:measure-monotonicity},
the whole functional does.

We can now verify that $R_i$ is a minimizer and prove the bounds \eqref{kp:eq:comparison-bounds}.
Weak continuity of $F_i$ on measures and
\eqref{kp:eq:measure-monotonicity} give
$F_i(R_i)\le F_i(\alpha)=\lim_\ell F_i(S_\ell)$.
Thus $R_i$ minimizes $F_i$ and, since $[T_i]\in\mathcal C_i$,
$F_i(R_i)\le\Mass(T_i)$. Together with \eqref{kp:eq:am-lower},
this proves the mass bounds and
\[
 \frac{1}{i^2}\sum_{j=1}^i
 \left(\int\varphi_j\,d(\|R_i\|-\|T_i\|)\right)^2
 =F_i(R_i)-\Mass(R_i)
 \le\Mass(T_i)-\Mass(R_i)\le\eps_i=\frac{1}{i^4}.
\]
This proves the test-integral bounds in
\eqref{kp:eq:comparison-bounds}; density of their span in $C(M)$
and the mass bound give $\|R_i\|\rightharpoonup\|V\|$.

It remains to prove the weighted stationarity.
Only $F_i$, not ordinary mass, is minimized on $\mathcal C_i$.
Set $w_i=w_i[\|R_i\|]$. Expanding the squares gives
\begin{equation}\label{kp:eq:auxiliary-expansion}
 F_i(S)-F_i(R_i)
 =\Mass_{w_i}(S)-\Mass_{w_i}(R_i)
 +\frac{1}{i^2}\sum_{j=1}^i
 \left(\int\varphi_j\,d(\|S\|-\|R_i\|)\right)^2,
\end{equation}
where $\Mass_{w_i}(S)=\int w_i\,d\|S\|$. Let $\Phi_t$ be a smooth
admissible flow supported in $\Omega$. Applying the homotopy formula
to an integral representative of $[R_i]$ gives
$\Flat([(\Phi_t)_\#R_i]-[R_i])\le C_i|t|$, since the flow sweeps
the boundary of that representative inside $\partial M$.
For either sign of sufficiently small $t$, this last flat step is
below $\delta_i$, and the mass is below $\Mass(T_i)+\delta_i$ by
continuity and $\Mass(R_i)\le\Mass(T_i)$. The endpoint therefore
belongs to $\mathcal C_i$. The last term in \eqref{kp:eq:auxiliary-expansion} is quadratic in time, so minimality gives zero first variation of $\Mass_{w_i}$, which is mass for $g_i=w_i^{2/k}g$.
\end{proof}

\subsection{The boundary force}
\label{kp:subsec:boundary-force}

We next control the boundary force and collar mass of the comparison currents.
In a fixed normal collar, let $d$ be inward distance and $\pi$ normal
projection, so $\nabla d=-\nu$. For a positive smooth weight $w$, put
$\operatorname{div}_{Q,w}X=\operatorname{div}_QX+
\langle\nabla\log w,X\rangle$; hence
$w\operatorname{div}_{Q,w}X=w\operatorname{div}_QX+
\langle\nabla w,X\rangle$.

\begin{lemma}
\label{kp:lem:boundary-force}
For each comparison current $R_i$, there is a nonnegative locally finite measure $\sigma_i$ on $\Omega\cap\partial M$ such that
\begin{equation}\label{kp:eq:boundary-force}
 \delta_{w_i}|R_i|(X):=
 \int_M\bigl(w_i\operatorname{div}_{T_xR_i}X
                  +\langle\nabla w_i,X\rangle\bigr)\,d\|R_i\|
 =\int_{\partial M}w_i\langle X,\nu\rangle\,d\sigma_i
\end{equation}
for every smooth vector field $X$ compactly supported in $\Omega$. The integral defining the left-hand side is also defined for fields pointing out of $M$; their flows need not preserve $M$.
For every nonnegative continuous function $b$ compactly supported in $\Omega\cap\partial M$,
\begin{equation}\label{kp:eq:weighted-collar-localization}
 \limsup_{\varepsilon\downarrow0}\frac1\varepsilon
 \int_{\{0<d<\varepsilon\}}(b\circ\pi)w_i\,d\|R_i\|
 \le\int_{\partial M}b w_i\,d\sigma_i.
\end{equation}
The same inequality holds with $b\circ\pi$ replaced by any continuous function $\widetilde b$ on the collar such that $\widetilde b|_{\partial M}=b$ and $\spt\widetilde b\subset\pi^{-1}(K)$ for some compact set $K\Subset\Omega\cap\partial M$.
\end{lemma}

\begin{proof}
Fix $i$ and write $R=R_i$, $w=w_i$, and $F=F_i$. Constants $C_i$ may depend on this fixed index and the cutoff functions.

\noindent\textbf{Step 1: The boundary force formula.}\enspace
By Lemma~\ref{kp:lem:comparison-minimizers}, $|R|$ is stationary with free boundary for the metric $g_i=w^{2/k}g$; that is, its first variation for $g_i$ vanishes on vector fields tangent to $\partial M$. De Masi's boundary first-variation formula \cite[Corollary 4.7]{DeMasi2021} describes the first variation of such a varifold on arbitrary vector fields. It requires a free-boundary first variation with an $L^1$ mean-curvature term, but assumes neither rectifiability nor integrality.

The formula is stated in Euclidean space, and we apply it locally as follows. Extend $(M,g_i)$ to a closed Riemannian manifold and embed the extension isometrically in Euclidean space by Nash's theorem, as in \cite[Section~2.2]{LZ2021}. Under the tubular projection onto this extension, the inverse image of $\partial M$ is an ambient hypersurface whose unit normal along $\partial M$ is the outward $g_i$-unit normal $\nu_i$. If an ambient vector field is tangent to this hypersurface, then its tangential projection onto $M$ is tangent to $\partial M$; thus intrinsic stationarity supplies the ambient first-variation identity required by the formula. The additional mean curvature due to the embedding is bounded and normal to $M$. It is therefore tangent to the extended boundary, as De Masi's formula requires, and it contributes zero against fields tangent to $M$.

The intrinsic mean curvature vanishes by weighted stationarity. The boundary-curvature term in De Masi's formula also vanishes upon integration, since $\|R\|_{g_i}(\partial M)=0$. Thus there is a nonnegative boundary measure $\widetilde\sigma_i$ such that
\[
 \delta_w|R|(X)
 =\int_{\partial M}\langle X,\nu_i\rangle_{g_i}\,d\widetilde\sigma_i
 =\int_{\partial M}w^{1/k}\langle X,\nu\rangle\,d\widetilde\sigma_i.
\]
Here $\|R\|_{g_i}=w\|R\|$ is the mass measure for $g_i$, and $\nu_i=w^{-1/k}\nu$. Defining
$d\sigma_i=w^{1/k-1}\,d\widetilde\sigma_i$
proves \eqref{kp:eq:boundary-force}. Write $\sigma=\sigma_i$ below.

\noindent\textbf{Step 2: The upper collar bound.}\enspace
Fix $0\le\varphi\in C_c^\infty(\Omega\cap\partial M)$ and use the same symbol for its extension $\varphi\circ\pi$. For small $t>0$, put
\[
 E_t=\{0<d\le t\varphi\},\quad
 m(t)=\int_{E_t}w\,d\|R\|,\quad
 B_\varphi=\int_{\partial M}\varphi w\,d\sigma.
\]
We will prove
\begin{equation}\label{kp:eq:collar-upper}
 \limsup_{t\downarrow0}\frac{m(t)}t\le B_\varphi.
\end{equation}
Choose a smooth cutoff $\chi$ of the signed normal coordinate, equal to one near zero and supported in a sufficiently small interval that $\chi(d)\varphi\nu$ has compact support in $\Omega$ when restricted to $M$. Let $\Psi_t$ be its ambient flow, and let $P_t$ be this flow stopped at its first boundary contact. Near $\partial M$,
\[
 \Psi_t(y,d)=(y,d-t\varphi(y)),\quad
 P_t(y,d)=(y,\max\{d-t\varphi(y),0\}).
\]
Thus $P_t$ fixes $\partial M$, maps $E_t$ into $\partial M$, and agrees with $\Psi_t$ on $M^\circ\setminus E_t$.

Set $S_t=((P_t)_\#R)_{\rm can}$. To estimate its relative flat distance from $R$, choose an integral representative $\widetilde R$ of $[R]$. Define $H(s,x)=P_s(x)$ for $0\le s\le t$ and let $Q_t=H_\#([0,t]\times\widetilde R)$ be the homotopy current. Since $H$ is constant in $s$ on $\partial M$ and $\spt\partial\widetilde R\subset\partial M$, the homotopy formula gives $\partial[Q_t]=[S_t]-[R]$. The time and spatial derivatives of $H$ are uniformly bounded for small $t$, so
$\Flat([S_t]-[R])\le\Mass(Q_t)\le Ct\Mass(\widetilde R)\le C_i t$.
The area formula on $M^\circ\setminus E_t$ also gives $\Mass(S_t)\le\Mass(R)+C_i t$. Since $\Mass(R)\le\Mass(T_i)$, for small $t$ this endpoint can be appended to a sequence reaching $[R]$. Hence $[S_t]\in\mathcal C_i$ and $F(S_t)\ge F(R)$.

Extend $w$ smoothly across the boundary for the following ambient area calculation. Equation~\eqref{kp:eq:boundary-force} gives
$\Mass_w((\Psi_t)_\#R)-\Mass_w(R)=tB_\varphi+O_i(t^2)$,
where this mass includes the portion outside $M$. The weighted Jacobian is $1+O_i(t)$, so subtracting the image of $E_t$ yields
\begin{equation}\label{kp:eq:stopped-area}
 \Mass_w(S_t)-\Mass_w(R)
 \le tB_\varphi-m(t)+C_i\bigl(t^2+tm(t)\bigr).
\end{equation}

Let $J_t(x)$ be the unweighted $k$-Jacobian of $\Psi_t$ on $T_xR$.
On $M^\circ\setminus E_t$, $P_t$ agrees with the injective map
$\Psi_t$ and has image in the interior; the rest is removed by
canonicalization. There is therefore no cancellation in the following
area formula. For each test function $\varphi_j$, it gives
\[
 \int\varphi_j\,d(\|S_t\|-\|R\|)
 =\int_{M\setminus E_t}
   [\varphi_j(\Psi_t(x))J_t(x)-\varphi_j(x)]\,d\|R\|(x)
   -\int_{E_t}\varphi_j\,d\|R\|.
\]
Using displacement $O(t)$, $J_t=1+O_i(t)$, and $w\ge1/2$ gives
$\left|\int\varphi_j\,d(\|S_t\|-\|R\|)\right|\le C_i\|\varphi_j\|_{C^1}(t+m(t))$.
Since $\frac{1}{i^2}\sum_{j=1}^i\|\varphi_j\|_{C^1}^2\le1$, substitution in \eqref{kp:eq:auxiliary-expansion}, together with \eqref{kp:eq:stopped-area}, gives
\[
 0\le F(S_t)-F(R)
 \le tB_\varphi-m(t)+C_i\bigl(t^2+tm(t)+m(t)^2\bigr).
\]
Since $m(t)\to0$, choose $t$ with $C_i(t+m(t))\le1/2$.
Then $m(t)/2\le tB_\varphi+C_it^2$, so $m(t)=O_i(t)$.
Substitution gives $m(t)\le tB_\varphi+O_i(t^2)$ and
\eqref{kp:eq:collar-upper}.

\noindent\textbf{Step 3: Constant-thickness localization.}\enspace
Define a boundary measure by
\[
 \alpha_\varepsilon(A)=\frac1\varepsilon
 \int_{\{0<d<\varepsilon,\ \pi(x)\in A\}}w\,d\|R\|
\]
for every Borel set $A\subset\Omega\cap\partial M$. Choosing $\varphi=1$ near any fixed compact set $A$ gives $\alpha_\varepsilon(A)\le m_\varphi(\varepsilon)/\varepsilon$, where the subscript records the choice of $\varphi$. By \eqref{kp:eq:collar-upper}, these measures have uniformly bounded mass on compact subsets. Every sequence $\varepsilon\downarrow0$ therefore has a subsequence converging weakly on compact subsets to a Radon measure $\alpha$.

Let $A\subset O\Subset\Omega\cap\partial M$, with $A$ compact and $O$ open. Choose smooth $0\le\xi,\varphi\le1$ such that $\xi=1$ on $A$, $\spt\xi\subset\{\varphi=1\}$, and $\spt\varphi\subset O$. Along the convergent subsequence,
\[
 \alpha(A)\le\int\xi\,d\alpha
 =\lim\int\xi\,d\alpha_\varepsilon
 \le\limsup\frac{m_\varphi(\varepsilon)}\varepsilon
 \le\int\varphi w\,d\sigma\le\int_O w\,d\sigma.
\]
Outer regularity of $w\sigma$ gives $\alpha(A)\le(w\sigma)(A)$
for compact $A$; inner regularity then gives $\alpha\le w\sigma$
on all Borel sets.

Now let $b\ge0$ be continuous with compact support. The quantities $\int b\,d\alpha_\varepsilon$ are bounded. Choose a sequence $\varepsilon\downarrow0$ realizing their upper limit and then a subsequence with weak limit $\alpha$. Weak convergence and the preceding measure inequality give
\[
 \limsup_{\varepsilon\downarrow0}\frac1\varepsilon
 \int_{\{0<d<\varepsilon\}}(b\circ\pi)w\,d\|R\|
 =\int b\,d\alpha\le\int b w\,d\sigma,
\]
which proves \eqref{kp:eq:weighted-collar-localization}.
Finally, if $\widetilde b$ is a continuous extension supported over a fixed compact set $A$, uniform continuity gives
$\sup_{0\le d\le\varepsilon,\,\pi(x)\in A}|\widetilde b(x)-b(\pi(x))|=o(1)$.
The difference between the two rescaled integrals is therefore bounded by $o(1)\alpha_\varepsilon(A)=o(1)$, proving the extension assertion.
\end{proof}

\subsection{Passage to the limit}
\label{kp:subsec:limit}

The boundary forces pass to the limit on the fixed set
$\Omega'\Subset\Omega$.

\begin{lemma}
\label{kp:lem:limit-form}
The varifold $V$ is integral in $\Omega'$. There is a nonnegative Radon measure $\sigma$ on $\Omega'$, supported on $\partial M\cap\Omega'$, such that
\begin{equation}\label{kp:eq:limit-force}
 \delta V(X)=\int\langle X,\nu\rangle\,d\sigma
\end{equation}
for every smooth vector field $X$ compactly supported in $\Omega'$.
\end{lemma}

\begin{proof}
We first show that the first variations of the varifolds $|R_i|$ are uniformly bounded.
Applying \eqref{kp:eq:boundary-force} to $X/w_i$ and expanding
the quotient derivative gives, with
$H_i=\pi_{(T_xR_i)^\perp}\nabla\log w_i$,
\begin{equation}\label{kp:eq:approximant-force}
 \delta|R_i|(X)
 =-\int\langle H_i,X\rangle\,d\|R_i\|
   +\int\langle X,\nu\rangle\,d\sigma_i.
\end{equation}
Here $\|H_i\|_\infty\le C/i$ by $w_i\ge1/2$ and the $C^1$
bound in Lemma~\ref{kp:lem:comparison-minimizers}.

Choose $0\le\chi\le1$ smooth, compactly supported in the normal collar inside $\Omega$, and equal to one near $\overline{\Omega'}\cap\partial M$. Testing \eqref{kp:eq:approximant-force} with $\chi\nu$ yields
\[
 \begin{aligned}
 \int\chi\,d\sigma_i
 &=\int\bigl(\operatorname{div}_{T_xR_i}(\chi\nu)
                  +\chi\langle H_i,\nu\rangle\bigr)\,d\|R_i\|\\
 &\le\int\bigl(|\nabla\chi|
                   +k\chi|\nabla\nu|+\chi|H_i|\bigr)\,d\|R_i\|\\
 &\le C_\chi\Mass(R_i).
 \end{aligned}
\]
Here $C_\chi$ is independent of $i$. Since $\chi=1$ near
$\overline{\Omega'}\cap\partial M$, this bounds
$\sigma_i(\Omega'\cap\partial M)$; \eqref{kp:eq:approximant-force}
then gives, for $X$ compactly supported in $\Omega'$,
\[
 \left|\delta|R_i|(X)\right|
 \le\left[\int_{\Omega'}|H_i|\,d\|R_i\|
       +\sigma_i(\Omega'\cap\partial M)\right]\|X\|_\infty
 \le C\Mass(R_i)\|X\|_\infty.
\]
Thus the full first variations are uniformly bounded on $\Omega'$.

We can now prove that $V$ is integral.
The following results are applied in local isometric embeddings of a smooth extension of $(M,g)$, as in Lemma~\ref{kp:lem:boundary-force}. The embedding curvature adds a bounded mean-curvature term and therefore preserves the local first-variation bounds.

Allard's integral compactness theorem \cite[Theorem 6.4]{Allard1972},
\cite[Chapter 8, Theorem 5.8 and Remark 5.9]{Simon2014} now gives
a subsequence $|R_i|\to W$ with $W$ integral on $\Omega'$.
Lemma~\ref{kp:lem:comparison-minimizers} gives $\|W\|=\|V\|$ there.

Apply De Masi's formula \cite[Corollary 4.7]{DeMasi2021} separately to the original $V$, using its free-boundary stationarity in $U$. The embedding argument above supplies its bounded ambient mean-curvature term: when several planes occur over a point, average the trace of the embedding's second fundamental form over those planes using disintegration of $V$ over $\|V\|$. The formula gives locally bounded full first variation, with the boundary-curvature term retained because $V$ may carry boundary mass. Since $W$ is integral, $\|V\|=\|W\|$ has positive finite $k$-density almost everywhere. The rectifiability theorem \cite[Chapter 8, Theorem 5.5]{Simon2014} therefore makes $V$ rectifiable.

Both $V$ and $W$ are now rectifiable with the same mass measure, and
the approximate tangent plane of a rectifiable varifold is determined
almost everywhere by its mass measure. The rectifiable-varifold
representation \cite[Chapter 4, Section 1, discussion after (1.2);
Chapter 8, Section 1, p.~206]{Simon2014} therefore gives $W=V$ on
$\Omega'$, so $V$ is integral there. Rectifiability is essential:
equality of mass measures alone would not identify arbitrary varifolds.

Finally, we pass to the limit in the boundary forces.
After a further subsequence, $\sigma_i\rightharpoonup\sigma\ge0$
locally on $\Omega'$, with $\spt\sigma\subset\partial M$
\cite[Chapter 1, Theorem 4.16]{Simon2014}.
Since $\|H_i\|_\infty\Mass(R_i)\to0$, the limit in
\eqref{kp:eq:approximant-force} is \eqref{kp:eq:limit-force}.
\end{proof}

\begin{proof}[Completion of the proof of Proposition~\ref{loc:prop:integral}]
Take an open exhaustion $U_j\Subset U_{j+1}$ of $U$.
Lemma~\ref{kp:lem:limit-form}, applied with $\Omega'=U_j$ and
$\Omega=U_{j+1}$, gives integrality in $U_j$ and a boundary measure
there. Testing with $X=\varphi\nu$ near $\partial M$, for
$\varphi\in C_c^\infty(U_j)$, shows that the measure on
$U_{j+1}\cap\partial M$ restricts to the one on $U_j\cap\partial M$.
Thus $V$ is integral in $U$, and these measures define a unique
nonnegative Radon measure $\sigma$ on $U\cap\partial M$. Every smooth
field compactly supported in $U$ has support in some $U_j$, so the
first-variation formula holds for each such field.
\end{proof}

\begin{corollary}
\label{loc:cor:integral}
Let $2\le k\le n-1$. Every $k$-varifold in $M$ that is stationary with free boundary and almost minimizing modulo two in small relative annuli is integral.
\end{corollary}

\begin{proof}
Let $V$ be such a varifold, fix $p\in M$, and let $r_p$ be as in Definition~\ref{def:almostminimizing}. By Proposition~\ref{loc:prop:integral}, $V$ is integral in every annulus $\An(p;s,r)$ with $0<s<r<r_p$. These annuli cover $M\cap B_{r_p}(p)\setminus\{p\}$. The monotonicity formula for stationary varifolds, \cite[Section 5]{Allard1972} at interior points and \cite[Theorem 2.3]{LZ2021} at boundary points, gives $\|V\|(B_r(p))\le Cr^k$ for all small $r>0$, and hence $\|V\|(\{p\})=0$. Thus $V$ is integral in $M\cap B_{r_p}(p)$. Finitely many of these balls cover $M$, so $V$ is integral.
\end{proof}

\section{Boundary trace and geometric consequences}
\label{sec:kflat}

We prove a local trace identity
(Proposition~\ref{kp:local-trace}) and then apply it on small annuli to
obtain Theorem~\ref{thm:kflat} and
Corollaries~\ref{cor:zeroboundary}--\ref{cor:dirichlet}.

The model computation is the following. For a smooth contact disk
$D\subset\partial M$ and $0\le\varphi\in C_c^\infty(D)$, the inward
velocity $-\varphi\nu$ fixes the edge and changes area at the rate
\[
 -\int_D\varphi\operatorname{tr}_{TD}\SFF_{\partial M}\,d\HH^k.
\]
Since $\nu\perp TD$, no first-order $\nabla\varphi$ term occurs, so
local mass minimality would exclude positive trace. In general,
however, the varifold is only a limit of almost minimizing cycles, and
an approximate tangent plane need not lie in an open smooth disk of the
support. We therefore work with the comparison currents $R_i$ of
Section~\ref{sec:local}, which are selected by auxiliary functionals and
preserve the prescribed limiting mass measure, including its boundary
part. By Proposition~\ref{loc:prop:integral}, the first variation of $V$
in $U$ has the form $\delta V=\nu\,\sigma$ for a nonnegative measure
$\sigma$ on $U\cap\partial M$. Using this first-variation formula and Menne's second-order rectifiability
theorem \cite[Theorem~4.8]{Menne2013}, we show in
Theorem~\ref{kp:thm:nonnegative-trace} that the boundary trace equals the
density of the absolutely continuous part of $\sigma$ with respect to
$\|V\|$. We then test the first variation and apply Menne's Sobolev
inequality \cite[Corollary~5.28]{Menne2016} to prove the height and tilt
estimates of Lemma~\ref{kp:lem:excess} by iteration. An inward deformation fixing the
boundary pointwise, combined with a boundary-preserving expansion, then
makes the area decrease dominate the collar cost, and
\eqref{kp:eq:finite-area-bounds} excludes positive trace. No prior
smoothness or local mass minimality of the limit is used.

The section is organized as follows. Subsection~\ref{kp:subsec:trace-sign} proves that the trace is nonnegative (Theorem~\ref{kp:thm:nonnegative-trace}), establishes height and tilt estimates at boundary points (Lemma~\ref{kp:lem:excess}), and excludes positive trace (Lemma~\ref{kp:lem:positive-trace}). Subsection~\ref{sec:consequences} derives the global statements, and Subsection~\ref{sec:examples} shows by examples that the hypotheses are needed. We use the notation and conventions of Section~\ref{sec:local} for $\Omega'\Subset\Omega\Subset U$, the comparison currents $R_i$, their weights $w_i$, and their boundary measures $\sigma_i$. On $\Omega'$ we use the restriction of $\sigma$ from Proposition~\ref{loc:prop:integral}.

By the tangent-plane and differentiation theorems
\cite[Chapter 3, Theorem 1.9]{Simon2014},
\cite[Sections 2.9 and 3.2]{Federer1969}, an integral varifold has a
unique integer-multiplicity planar tangent almost everywhere;
on its boundary part this plane lies in $T\partial M$ almost everywhere.

\begin{proposition}
\label{kp:local-trace}
Let $2\le k\le n-1$, let $U\subset M$ be relatively open, and let $V$
be a $k$-varifold in $M$ that is almost minimizing modulo two in $U$.
Then
\begin{equation}\label{kp:eq:local-trace}
 \operatorname{tr}_{T_xV}\SFF_{\partial M}=0
 \quad\text{for }\|V\|\text{-almost every }x\in U\cap\partial M.
\end{equation}
The measure $\sigma$ of Proposition~\ref{loc:prop:integral} is singular
with respect to $\|V\|$ on $U$.
\end{proposition}

Recall that $V$ is integral in $U$ by Proposition~\ref{loc:prop:integral}, so the tangent planes in \eqref{kp:eq:local-trace} exist almost everywhere.

\subsection{Vanishing of the boundary trace}
\label{kp:subsec:trace-sign}
\label{kp:subsec:contact}
\label{kp:subsec:finite-comparison}

We first prove nonnegativity and construct $C^2$ boundary patches;
contact estimates will then exclude positive trace.

\begin{theorem}
\label{kp:thm:nonnegative-trace}
Let $V$ and $\sigma$ be as in Lemma~\ref{kp:lem:limit-form}, and write $\sigma=q\,\|V\|\res\partial M+\sigma_s$ with $q\ge0$ and $\sigma_s\perp\|V\|$. Then
\[\operatorname{tr}_{T_xV}\SFF_{\partial M}=q(x)\ge0\quad\text{for }\|V\|\text{-almost every }x\in\Omega'\cap\partial M.\]
\end{theorem}

\begin{proof}
The decomposition of $\sigma$ into its absolutely continuous and singular parts exists by \cite[Chapter 1, Lemma 3.19 and Theorem 3.24(ii)]{Simon2014}, and $q\ge0$ because $\sigma$ is nonnegative. We have to identify the trace with the density $q$. To this end we show that, at almost every boundary point, $V$ is tangent to second order to a $C^2$ submanifold of $\partial M$; the mean curvature of such a submanifold in the direction $\nu$ is minus the trace.

We begin with the generalized mean curvature of $V$.
Define
\[H_V(x)=\begin{cases} 0,&x\in\Omega'\cap\Int M,\\ -q(x)\nu(x),&x\in\Omega'\cap\partial M. \end{cases}\]
Then \eqref{kp:eq:limit-force} gives, for every smooth field $X$ compactly supported in $\Omega'$,
\[\delta V(X)=-\int_{\Omega'}\langle H_V,X\rangle\,d\|V\| +\int_{\Omega'\cap\partial M}\langle X,\nu\rangle\,d\sigma_s.\]
This identifies the generalized mean curvature $H_V$, with
$\int_K|H_V|\,d\|V\|\le\sigma(K)<\infty$ for compact
$K\subset\Omega'$.

Next, apply \cite[Theorem 4.8]{Menne2013} in the local isometric embeddings from Lemma~\ref{kp:lem:limit-form}. Its hypotheses are integrality and locally bounded full first variation; a singular part such as $\sigma_s$ is allowed. Denote the resulting countable family of $C^2$ embedded $k$-dimensional submanifolds by $N_j$. These are submanifolds of the embedding space, and
\[\|V\|\Bigl(\Omega'\setminus\bigcup_jN_j\Bigr)=0,\quad T_xN_j=T_xV,\quad H_{N_j}^{\mathrm E}(x)=H_V^{\mathrm E}(x),\]
where the last two identities hold at $\|V\|$-almost every $x\in N_j\cap\Omega'$. The superscript $\mathrm E$ denotes Euclidean mean curvature.

We now replace the submanifolds $N_j$, along the boundary part of $V$, by $C^2$ patches contained in $\partial M$.
Cover the $N_j$ by countably many $C^2$ embeddings
$\psi:O\subset\mathbb R^k\to\mathbb R^L$. In a tubular
neighborhood of the embedded boundary, let $\Pi$ be nearest-point
projection and put $E=\psi^{-1}(\partial M)$,
$\widehat\psi=\Pi\circ\psi$ and $F=\psi-\widehat\psi$.
At a Lebesgue density point $y\in E$, every unit direction $v$
is a limit of directions $(z_j-y)/|z_j-y|$ with $z_j\in E$ and
$z_j\to y$. Since $F|_E=0$, successive first- and second-order
Taylor expansions give $DF(y)v=0$ and $D^2F(y)[v,v]=0$.
Polarization therefore yields
\[
 \widehat\psi(y)=\psi(y),\quad
 D\widehat\psi(y)=D\psi(y),\quad
 D^2\widehat\psi(y)=D^2\psi(y).
\]
In particular, $D\widehat\psi(y)$ has rank $k$. The inverse
function theorem, applied to a nonsingular coordinate minor, makes
$\widehat\psi$ a $C^2$ graph embedding near $y$. Restricting to
rational balls gives countably many patches $S_\ell\subset\partial M$.
These patches cover the boundary part only up to measure zero;
no open portion of the support is asserted to be smooth.
The density theorem and local Lipschitz continuity of $\psi$ give
\[
 \|V\|\Bigl((\Omega'\cap\partial M)
                  \setminus\bigcup_\ell S_\ell\Bigr)=0.
\]
For each intersection $N_j\cap S_\ell$, repeat the same Taylor
argument with the $C^2$ graph retraction onto $S_\ell$ in place of
$\Pi$. Tangent planes and Euclidean mean curvatures agree almost
everywhere on the intersection. Since the $N_j$ cover $V$ up to
zero mass,
\[
 T_xS_\ell=T_xV,\quad H_{S_\ell}^{\mathrm E}(x)=H_V^{\mathrm E}(x)
 \quad\text{for }\|V\|\text{-almost every }x\in S_\ell\cap\Omega'.
\]

It remains to compare the trace with the density $q$.
The common tangent plane contributes the same trace of the
embedding's second fundamental form to both Euclidean mean
curvatures. Subtracting it gives $H_{S_\ell}=H_V=-q\nu$ in $M$.
Write $S=S_\ell$ and choose an orthonormal tangent frame $e_1,\ldots,e_k$.
Differentiating $\langle e_a,\nu\rangle=0$ gives
\[\langle H_S,\nu\rangle =\sum_{a=1}^k\langle\nabla_{e_a}e_a,\nu\rangle =-\sum_{a=1}^k\langle e_a,\nabla_{e_a}\nu\rangle =-\sum_{a=1}^k\SFF_{\partial M}(e_a,e_a),\]
where the last equality uses \eqref{eq:boundaryform}. Therefore
\begin{equation}\label{kp:eq:contact-curvature}
 h(x):=\operatorname{tr}_{T_xV}\SFF_{\partial M}
 =-\langle H_S,\nu\rangle=q(x)\ge0.
\end{equation}
The countable covering proves the assertion and gives $H_S=-h\nu$
on each common part.
\end{proof}

The next height and tilt estimates near a boundary patch $S$ control
the comparison errors in Lemma~\ref{kp:lem:positive-trace}.

In a tubular neighborhood of a smaller patch of $S$, let
$\pi_S$ be nearest-point projection and $u=\dist(\cdot,S)$.
Write $P(x)$ for projection onto $T_xV$ and $Q(x)$ for projection
onto $T_{\pi_S(x)}S$ transported along the short normal geodesic
to $x$. We use the Hilbert--Schmidt norm
$|P-Q|^2=\sum_{b=1}^n|(P-Q)v_b|^2$ for an orthonormal basis
$(v_b)$ of $T_xM$. Distances and geodesics lie in the smooth extension.

\begin{lemma}
\label{kp:lem:excess}
For $\|V\|$-almost every $p\in S$,
\begin{equation}\label{kp:eq:excess-conclusion}
 \int_{B_r(p)}u^2\,d\|V\|=o(r^{k+4}),\quad
 \int_{B_r(p)}|P(x)-Q(x)|^2\,d\|V\|=o(r^{k+2}).
\end{equation}
Here $B_r(p)$ is the ball of radius $r$ centered at $p$.
\end{lemma}

\begin{proof}
The proof is a Caccioppoli-type argument. Testing the first variation with the normal displacement field of $S$ bounds the tilt by the height. A Sobolev inequality, applied on the small set where $V$ leaves $S$, bounds the height by the tilt with a small factor. An iteration over dyadic radii then gives the decay.

We begin with two inequalities for test functions.
For a nonnegative compactly supported $C^1$ function $f$, put
$\nabla_Vf=P\nabla f$. Testing \eqref{kp:eq:limit-force} with
$f\nu$, using the fixed normal extension of $\nu$, gives
\begin{equation}\label{kp:eq:boundary-integral}
 \begin{aligned}
 \int f\,d\sigma&=\delta V(f\nu)\\
 &=\int\left[\langle\nabla_Vf,\nu\rangle+f\sum_{a=1}^k\langle\nabla_{e_a}\nu,e_a\rangle\right]d\|V\|\\
 &\le C\int(|\nabla_Vf|+f)\,d\|V\|.
 \end{aligned}
\end{equation}
Here $e_1,\ldots,e_k$ is an orthonormal basis of $T_xV$. Constants $C$ depend on the fixed neighborhood and may increase below, but are independent of $r$.
All $C^1$ test fields used below are allowed by smooth approximation:
their derivatives converge uniformly on a fixed compact set, and both
$\|V\|$ and $\sigma$ are finite there.

Let $B$ be the second fundamental form of the local isometric
embedding. Decomposing a Euclidean test field into tangential and
normal parts, \eqref{kp:eq:limit-force} gives
\[\delta_{\mathrm E}V(X)=\int\langle X,\nu\rangle\,d\sigma-\int\left\langle X,\sum_{a=1}^kB(e_a,e_a)\right\rangle d\|V\|.\]
Thus its total variation measure satisfies $\|\delta_{\mathrm E}V\|\le\sigma+C\|V\|$, because $B$ is bounded. Set $\beta=k/(k-1)$. In the Euclidean chart, integrality gives density at least one almost everywhere, and $f$ belongs to the zero-boundary Sobolev class in \cite[Corollary 5.28]{Menne2016}: it has compact support and admits smooth approximation in the relevant first-derivative norm. That result uses the full first-variation measure, so no absolute continuity of $\sigma$ is needed. It gives
\begin{equation}\label{kp:eq:contact-sobolev}
 \begin{aligned}
 \|f\|_{L^\beta(\|V\|)}
 &\le C\left(\int|\nabla_Vf|\,d\|V\|+\int f\,d\|\delta_{\mathrm E}V\|\right)\\
 &\le C\left(\int(|\nabla_Vf|+f)\,d\|V\|+\int f\,d\sigma\right)\\
 &\le C\int(|\nabla_Vf|+f)\,d\|V\|.
 \end{aligned}
\end{equation}

Next we bound the tilt by the height.
By measure differentiation and rectifiability, almost every $p\in S$ satisfies $\|V\|(B_r(p)\setminus S)=o(r^k)$. Fix such a point, write $B_r=B_r(p)$, and set
\[A(r)=\|V\|(B_r\setminus S),\quad H(r)=\int_{B_r}u^2\,d\|V\|,\quad T(r)=\int_{B_r}|P-Q|^2\,d\|V\|.\]
In particular, $A(r)=o(r^k)$. Let $w=\nabla(u^2/2)$. This is the normal displacement vector from $S$, so $|w|=u$ and $Qw=0$. Differentiating in normal coordinates gives $\nabla w=\id-Q+O(u)$, and therefore
\[\begin{aligned}|Pw|&\le u|P-Q|,\\ \operatorname{div}_{T_xV}w&=\operatorname{tr}\bigl(P(\id-Q)\bigr)+O(u)=\tfrac12|P-Q|^2+O(u).\end{aligned}\]
The last equality uses $\operatorname{rank}P=\operatorname{rank}Q=k$. Let $\rho$ be signed distance to $\partial M$ in a fixed collar, positive outside $M$. Thus $\rho=0$ and $\nabla\rho=\nu$ on $\partial M$.
For $x\in\partial M$ close to $S$, the normal geodesic
$\gamma(t)=\exp_x(-tw(x))$, $0\le t\le1$, joins $x$ to
$\pi_S(x)\in\partial M$ with speed $u(x)$. Taylor's formula gives
\[
 \begin{aligned}
 0&=\rho(\gamma(1))-\rho(\gamma(0))\\
  &=-\langle\nu(x),w(x)\rangle+\int_0^1(1-t)(\nabla^2\rho)_{\gamma(t)}(\dot\gamma(t),\dot\gamma(t))\,dt.
 \end{aligned}
\]
Hence $|\langle w,\nu\rangle|\le Cu^2$ on $\partial M$.

Choose $0\le\zeta\le1$, equal to one on $B_r$ and supported in $B_{2r}$, with $|\nabla\zeta|\le C/r$. We take $r<1$ and small enough that $B_{2r}$ lies in the fixed neighborhood. Put
\[f=\zeta^2u^2,\quad J=\int f\,d\|V\|,\quad I=\int\zeta^2|P-Q|^2\,d\|V\|.\]
The product rule and Cauchy--Schwarz give
\[
 \begin{aligned}
 \int|\nabla_Vf|\,d\|V\|
 &\le2\int\zeta u^2|\nabla\zeta|\,d\|V\|+2\int\zeta^2|Pw|\,d\|V\|\\
 &\le Cr^{-1}J^{1/2}H(2r)^{1/2}+2J^{1/2}I^{1/2}.
 \end{aligned}
\]
We use this bound twice: first to control the boundary term, and then in \eqref{kp:eq:contact-sobolev}.

Test \eqref{kp:eq:limit-force} with $\zeta^2w$. Using
$\operatorname{div}_{T_xV}w\ge\tfrac12|P-Q|^2-Cu$,
$|\langle w,\nu\rangle|\le Cu^2$, and
$|Pw|\le u|P-Q|$, Cauchy--Schwarz gives the first line below.
The gradient bound for $f$, $J\le H(2r)$, and
\eqref{kp:eq:boundary-integral} give the second:
\[
 \begin{aligned}
 \tfrac12 I
 &\le C\int_{B_{2r}}u\,d\|V\|+Cr^{-1}\sqrt{IH(2r)}+C\int f\,d\sigma\\
 &\le C\int_{B_{2r}}u\,d\|V\|+C(r^{-1}+1)\sqrt{IH(2r)}+C(r^{-1}+1)H(2r)\\
 &\le C\int_{B_{2r}}u\,d\|V\|+\tfrac14I+C^2(r^{-1}+1)^2H(2r)+C(r^{-1}+1)H(2r)\\
 &\le\tfrac14I+C\int_{B_{2r}}u\,d\|V\|+Cr^{-2}H(2r).
 \end{aligned}
\]
The penultimate line uses $2ab\le a^2+b^2$, and the last uses $r<1$. Since $\zeta=1$ on $B_r$, absorption yields
\begin{equation}\label{kp:eq:caccioppoli}
 T(r)\le I\le Cr^{-2}H(2r)+C\int_{B_{2r}}u\,d\|V\|.
\end{equation}

Conversely, the height is controlled by the tilt.
Since $f$ vanishes on $S$ and outside $B_{2r}$, H\"older's inequality and \eqref{kp:eq:contact-sobolev} give
\[
 \begin{aligned}J&\le A(2r)^{1/k}\|f\|_{L^\beta(\|V\|)}\\
 &\le CA(2r)^{1/k}\left[J^{1/2}\bigl(r^{-1}H(2r)^{1/2}+I^{1/2}\bigr)+J\right].\end{aligned}
\]
For small $r$, absorb the last term into the left side. Dividing by $J^{1/2}$ and squaring gives $J\le CA(2r)^{2/k}[r^{-2}H(2r)+I]$; if $J=0$, this is immediate. Also, $u=0$ on $S$, so Cauchy--Schwarz gives
\[\int_{B_{2r}}u\,d\|V\|\le\sqrt{A(2r)H(2r)}.\]
Combining these estimates with $H(r)\le J$ and \eqref{kp:eq:caccioppoli}, we obtain
\begin{equation}\label{kp:eq:small-support}
 H(r)\le CA(2r)^{2/k}
 \left[r^{-2}H(2r)+\sqrt{A(2r)H(2r)}\right].
\end{equation}
Young's inequality, in the form $r^2\sqrt{A(2r)H(2r)}\le[H(2r)+r^4A(2r)]/2$, now gives
\[\begin{aligned}H(r)&\le C\frac{A(2r)^{2/k}}{r^2}\left[H(2r)+r^4A(2r)\right]\\&=o(1)H(2r)+o(r^{k+4}).\end{aligned}\]

It remains to iterate this inequality.
Fix $\varepsilon>0$ and choose $r_0>0$ so that
$H(r)\le2^{-(k+5)}H(2r)+\varepsilon r^{k+4}$ for $0<r<r_0$.
Let $N=N(r)$ be the integer with $r_0/2<2^Nr\le r_0$. Dividing by $r^{k+4}$ and iterating gives
\[\frac{H(r)}{r^{k+4}}\le2^{-N}\frac{H(2^Nr)}{(2^Nr)^{k+4}}+2\varepsilon.\]
The first term tends to zero because $H(r_0)<\infty$ and $N\to\infty$. Since $\varepsilon$ is arbitrary, $H(r)=o(r^{k+4})$. Finally, \eqref{kp:eq:caccioppoli} gives
\[T(r)\le Cr^{-2}H(2r)+C\sqrt{A(2r)H(2r)}=o(r^{k+2}),\]
which proves \eqref{kp:eq:excess-conclusion}.
\end{proof}

We can now exclude positive trace by comparing the inward area
decrease with the cost of fixing the boundary.

\begin{lemma}
\label{kp:lem:positive-trace}
Let $V$ satisfy the hypotheses of Proposition~\ref{kp:local-trace}, and let $\Omega'\Subset\Omega\Subset U$ be as in Section~\ref{sec:local}. Then
\[
 \|V\|\bigl(\{x\in\Omega'\cap\partial M:h(x)>0\}\bigr)=0,
 \quad h(x)=\operatorname{tr}_{T_xV}\SFF_{\partial M}.
\]
\end{lemma}

\begin{proof}
Suppose otherwise. Choose a patch $S$ from the proof of Theorem~\ref{kp:thm:nonnegative-trace} on which this set has positive mass. By Lemma~\ref{kp:lem:excess} and measure differentiation applied to \eqref{kp:eq:contact-curvature}, we may choose $p\in S$ with $h(p)>0$ such that \eqref{kp:eq:excess-conclusion} holds and
\begin{equation}\label{kp:eq:chosen-contact-point}
 \|V\|(B_r)=\theta\omega_k r^k+o(r^k),\quad
 \sigma(B_r)=h(p)\theta\omega_k r^k+o(r^k),\quad \theta>0.
\end{equation}
Here $B_r=B_r(p)$ is the metric ball in the fixed smooth extension, intersected with $M$, and $\omega_k$ is the volume of the Euclidean unit $k$-ball. The second identity follows from $\sigma=h\|V\|\res\partial M+\sigma_s$: at almost every such $p$, both $\sigma_s(B_r)$ and $\sigma(B_r\setminus S)$ are $o(r^k)$, while $h|_S=\operatorname{tr}_{TS}\SFF_{\partial M}$ is continuous.

Before giving the details, we explain the idea. We compare the minimizers $R_i$ of Lemma~\ref{kp:lem:comparison-minimizers} with the currents
$S_i(s,t,\varepsilon)=(f_s\circ F_{t,\varepsilon})_\#R_i$,
which are obtained from two maps supported in $B_r(p)$. Let $G_t$
denote the uncut inward motion, of size $O(tr^2)$, and let
$F_{t,\varepsilon}$ be its cutoff in a boundary collar of width
$\varepsilon$. The map $f_s$ preserves $\partial M$ and, near $p$ in
coordinates adapted to the patch, dilates by $\lambda=(1+ks)^{1/k}$.
Write $B_i$ for $r^2$ times the localized boundary force of $R_i$.
Since $h(p)>0$, this quantity is of order $r^{k+2}$. For the uncut
motion $f_s\circ G_t\circ f_s^{-1}$ on the dilated current, the
leading decrease in weighted mass is $t\lambda^{k+1}B_i$. Replacing
$G_t$ by $F_{t,\varepsilon}$ fixes the boundary pointwise before
dilation and, by Lemma~\ref{kp:lem:boundary-force}, costs at most
$t\lambda^kB_i$ to leading order. The difference between these leading
terms is $t\lambda^k(\lambda-1)B_i\ge c\,st\,r^{k+2}$. The dilation itself is almost free: its Jacobian is affine in $s$ along the model patch, and $R_i$ is close to this patch by Lemma~\ref{kp:lem:excess}, so it increases the weighted mass by at most $o(1)\,s^2r^{k+2}$ as $r\to0$. Choosing $t=\eta s$ with $\eta$ and $r$ small, we find that the net gain dominates all errors. The resulting current lies in the comparison class $\mathcal C_i$ and has a smaller value of $F_i$ than $R_i$, which is the desired contradiction. Since the argument is long, we divide it into six steps.

\noindent\textbf{Step 1: Coordinates and supports.}\enspace
Flattening the $C^2$ patch $S$ would give only $C^1$ metric
coefficients. For the Jacobian estimates, we instead use a smooth
$S_0\subset\partial M$ agreeing with $S$ through second order at $p$.
Choose a smooth boundary parametrization $\kappa:W\to\partial M$, where $W$ is a neighborhood of $0$ in $\mathbb R^{n-1}$, such that $\kappa(0)=p$ and $D\kappa(0)(\mathbb R^k\times\{0\})=T_pS$. After restricting to a neighborhood of $p$, write
$S=\{\kappa(\xi,f(\xi)):\xi\in B^k_\rho(0)\}$, where $B^k_\rho(0)=\{\xi\in\mathbb R^k:|\xi|<\rho\}$ and $f:B^k_\rho(0)\to\mathbb R^{n-1-k}$ is $C^2$, with $f(0)=0$ and $Df(0)=0$. Define
\begin{equation}\label{kp:eq:quadratic-model}
 f_2(\xi)=\tfrac12D^2f(0)[\xi,\xi],\quad
 S_0=\{\kappa(\xi,f_2(\xi)):\xi\in B^k_\rho(0)\}.
\end{equation}
Shrink $\rho$ so that the smooth graph $S_0$ lies in the chart.
Taylor's theorem gives $f(\xi)-f_2(\xi)=o(|\xi|^2)$ and $Df(\xi)-Df_2(\xi)=o(|\xi|)$.
When $k=n-1$, the graph functions take values in $\mathbb R^0=\{0\}$, so $S_0$ is simply the corresponding boundary patch.
Now choose geodesic coordinates $v$ on $S_0$ centered at $p$, normal coordinates $z$ to $S_0$ within $\partial M$, and inward normal distance $d$. For each $i$, put
\[
 b_i(v)=w_i(v)\sqrt{\det g_{S_0}(v)},\quad
 y_1=\int_0^{v_1}\frac{b_i(\tau,v_2,\ldots,v_k)}{b_i(0)}\,d\tau,
 \quad y_j=v_j\ (j\ge2).
\]
The coordinates $x=(y,z,d)$ satisfy $x(p)=0$ and $w_i\,d\HH^k|_{S_0}=w_i(p)\,dy$. Their differentials at $p$ are orthonormal. Since $w_i\to1$ smoothly, the coordinate maps and their inverses converge smoothly. Shrink their common domain so that
$\tfrac12|x|\le\dist(p,x)\le2|x|$,
and choose $0<r<1$ such that $\overline B_r\subset\Omega'$ lies in this coordinate domain.
Let $L_{z,d}$ denote the submanifold obtained by varying $y$ with $z,d$ fixed. For a $k$-plane $Q\subset T_xM$, set
$u(x)=\dist(x,S_0)$ and $a(x,Q)=|\pi_Q-\pi_{T_xL_{z,d}}|$.
The submanifolds $L_{z,d}$ and the functions $u,a$ are independent of the choice of $y$. Replacing $S$ by $S_0$ changes distances by $o(r^2)$ and tangent projections by $o(r)$ on $B_r$. The plane $T_xL_{z,d}$ differs by $O(u(x))$ from $TS_0$ transported along the normal geodesic to $x$. Thus \eqref{kp:eq:excess-conclusion} implies
\begin{equation}\label{kp:eq:finite-contact-errors}
 \int_{B_r}u^2\,d\|V\|=o(r^{k+4}),\quad
 \int_{B_r}a(x,T_xV)^2\,d\|V\|=o(r^{k+2}).
\end{equation}
Choose $\zeta_r\in C_c^\infty(B_r)$ with $0\le\zeta_r\le1$ and $\zeta_r=1$ on $B_{r/2}$. For every continuous function $\varphi(x,Q)$ on $G_k(\overline B_r)$ and every continuous function $\psi$ on $\overline B_r\cap\partial M$, Lemma~\ref{kp:lem:limit-form} gives
\begin{equation}\label{kp:eq:cutoff-convergence}
 \begin{split}
 \lim_{i\to\infty}\int_{B_r}\zeta_r(x)\varphi(x,T_xR_i)\,d\|R_i\|(x)
 &=\int_{B_r}\zeta_r(x)\varphi(x,T_xV)\,d\|V\|(x),\\
 \lim_{i\to\infty}\int_{B_r\cap\partial M}\zeta_r\psi\,d\sigma_i
 &=\int_{B_r\cap\partial M}\zeta_r\psi\,d\sigma.
 \end{split}
\end{equation}
Taking $\varphi=\psi=1$ and using \eqref{kp:eq:chosen-contact-point}, we obtain
\[
 \limsup_{i\to\infty}\bigl[\|R_i\|(B_{r/2})+\sigma_i(B_{r/2})\bigr]
 \le\int\zeta_r\,d\|V\|+\int\zeta_r\,d\sigma
 \le\|V\|(B_r)+\sigma(B_r)\le Cr^k.
\]
Hence $\|R_i\|(B_{r/2})+\sigma_i(B_{r/2})\le Cr^k$ for fixed small $r$ and sufficiently large $i$, after enlarging $C$. Formula~\eqref{kp:eq:cutoff-convergence} also holds for tests $\varphi_i,\psi_i$ converging uniformly to $\varphi,\psi$: for example,
\[
 \left|\int\zeta_r(\varphi_i-\varphi)\,d|R_i|\right|
 \le\sup_{G_k(\overline B_r)}|\varphi_i-\varphi|
       \int\zeta_r\,d\|R_i\|\longrightarrow0.
\]
The boundary integral is treated by the same inequality with $\psi_i-\psi$ and $\sigma_i$.
Constants $C$ below depend on the fixed charts and cutoff functions, but not on $i,r,s,t,\varepsilon$; $C_r$ may depend on $r$.

\noindent\textbf{Step 2: The comparison maps.}\enspace
We define the dilation $f_s$, the inward maps $G_t$ and $F_{t,\varepsilon}$, and the quantities $A_i(s,t)$ and $B_i$ in terms of which the estimates are expressed.
Fix a smooth function $\chi:[0,\infty)\to[0,1]$ with
$\chi=1$ on $[0,1]$ and $\chi=0$ on $[2,\infty)$.
For $R\ge0$, define
\[
 \chi_r(R)=\chi(8R/r),\quad
 \alpha_s(R)=(1+ks\chi_r(R))^{1/k},\quad
 f_s(x)=\alpha_s(|x|)x,
\]
and put $\lambda=(1+ks)^{1/k}$ and $\rho_r(x)=\chi(16|x|/r)$.
Choose $s_0\in(0,1)$ so that $\alpha_s(R)+R\alpha_s'(R)>0$ and
$s/2\le\lambda-1\le s$ for $0\le s\le s_0$.
Then $f_s$ is a diffeomorphism preserving $\partial M$, and
$f_s(x)=\lambda x$ on $\spt\rho_r$.

For $0\le t\le1$, put
\[
 Y=r^2\rho_r\partial_d=-r^2\rho_r\nu,\quad
 G_t(y,z,d)=(y,z,d+tr^2\rho_r(y,z,d)).
\]
To fix $\partial M$ pointwise, define
\begin{equation}\label{kp:eq:inward-map-comparison}
 \begin{split}
 q_\varepsilon(d)&=\min\{d/\varepsilon,1\},\\
 F_{t,\varepsilon}(y,z,d)
 &=(y,z,d+tr^2\rho_r(y,z,d)q_\varepsilon(d)),\\
 S_i(s,t,\varepsilon)&=(f_s\circ F_{t,\varepsilon})_\#R_i.
 \end{split}
\end{equation}
Let $G_t^d$ and $F_{t,\varepsilon}^d$ denote the last coordinates
of these maps. For $d\ne\varepsilon$,
\[
 \begin{split}
 \partial_dG_t^d&=1+tr^2\partial_d\rho_r,\\
 \partial_dF_{t,\varepsilon}^d
 &=1+tr^2q_\varepsilon\partial_d\rho_r
      +tr^2\rho_rq_\varepsilon'
 \ge1-16r\|\chi'\|_\infty.
 \end{split}
\]
Reduce $r$ so that $16r\|\chi'\|_\infty\le1/2$.
Both derivatives are then positive. Since $y,z$ are unchanged,
the maps are injective. They fix every point outside
$\spt\rho_r$, so $G_t(\spt\rho_r)\subset\spt\rho_r$.
Moreover, $F_{t,\varepsilon}$ is a bi-Lipschitz homeomorphism
of $M$ fixing $\partial M$ pointwise. Thus
$S_i(s,t,\varepsilon)$ is canonical
(that is, $\|S_i(s,t,\varepsilon)\|(\partial M)=0$),
and its full mass is given by the area formula.

For a linear map $A$ on a $k$-plane $Q$ with orthonormal basis
$(e_1,\ldots,e_k)$, write
$J_k(A|Q)=\sqrt{\det(\langle Ae_a,Ae_b\rangle)_{a,b=1}^k}$. Set
\[
 A_i(s,t)=\Mass_{w_i}((f_s\circ G_t)_\#R_i),\quad
 B_i=r^2\int w_i\rho_r\,d\sigma_i.
\]
Since $\rho_r=1$ on $B_{r/32}(p)$ and
$\spt\rho_r\subset\overline{B_{r/4}(p)}$
(these are metric balls; the coordinate norm satisfies
$\tfrac12|x|\le\dist(p,x)\le2|x|$),
\eqref{kp:eq:chosen-contact-point} and
\eqref{kp:eq:cutoff-convergence} give
$cr^{k+2}\le B_i\le Cr^{k+2}$;
here the coordinate cutoffs converge smoothly for fixed $r$.

\noindent\textbf{Step 3: Expansion.}\enspace
We show that the dilation alone changes the weighted mass only to second order in $s$, and with a small constant; the result is \eqref{kp:eq:expansion-bound}.
Write the coordinates as $x=(x^1,\ldots,x^n)$ and define
\[
 X(x)=\chi_r(|x|)\sum_{a=1}^n
 x^a\left.\frac{\partial}{\partial x^a}\right|_x,
 \quad
 j_{i,s}(x,Q)
 =\frac{w_i(f_s(x))}{w_i(x)}J_k(Df_s(x)|Q).
\]
Here $\partial/\partial x^a$ are the coordinate tangent fields
on $M$, and $X$ is extended by zero outside the coordinate
neighborhood. Since
$\left.\partial_sf_s(x)\right|_{s=0}=X(x)$
and $X$ is tangent to $\partial M$,
\eqref{kp:eq:boundary-force} gives
$\left.\partial_sA_i(s,0)\right|_{s=0}=\delta_{w_i}|R_i|(X)=0$.

On $S_0=\{z=d=0\}$, the identity
$w_i\,d\HH^k|_{S_0}=w_i(p)\,dy$
means that $j_{i,s}(x,T_xS_0)$ equals the Euclidean Jacobian
of $y\mapsto\alpha_s(|y|)y$. Write $R=|y|=|x|$ and let a prime denote
differentiation with respect to $R$. The radial derivative of this map
is $\alpha_s+R\alpha_s'$, while its derivatives in directions
perpendicular to $y$ are multiplication by $\alpha_s$. Since
$\alpha_s^k=1+ks\chi_r(R)$ and hence
$\alpha_s^{k-1}\alpha_s'=s\chi_r'(R)$, we obtain
\begin{equation}\label{kp:eq:model-jacobian}
 j_{i,s}(x,T_xS_0)
 =\alpha_s^{k-1}(\alpha_s+R\alpha_s')
 =1+s\bigl(k\chi_r(R)+R\chi_r'(R)\bigr);
\end{equation}
the formula at $R=0$ follows by continuity.
This exact affine dependence on $s$ removes the quadratic expansion
cost on the model patch. The remaining cost must therefore be
estimated by height and tilt errors, rather than by the total mass.

Write $\xi=(z,d)$ and $x_0=(y,0,0)$.
For an $(n-k)\times k$ matrix $\ell$, with
$|\ell|^2=\sum_{\beta,a}\ell_{\beta a}^2$, define
\[
 Q_\ell(x)=\operatorname{span}\left\{
 \partial_{y_a}
 +\sum_{\beta=1}^{n-k}\ell_{\beta a}\partial_{\xi_\beta}
 :1\le a\le k\right\}.
\]
At $x_0$, $Df_s(x_0)$ preserves the orthogonal tangent
and normal subspaces. Thus both Gram matrices contain
only constant and quadratic terms in $\ell$.
The weighted Jacobian is even in $\ell$, so its first
derivative with respect to the plane $Q$ vanishes:
\begin{equation}\label{kp:eq:plane-first-derivative}
 \left.D_\ell j_{i,s}(x_0,Q_\ell(x_0))\right|_{\ell=0}=0.
\end{equation}

Fix $x=(y,\xi)$ and a plane $Q=Q_\ell(x)$ with
$a(x,Q)$ small. The coordinate and metric norms are
uniformly equivalent, so $|\xi|\le Cu(x)$ and
$|\ell|\le Ca(x,Q)$. For $0\le\tau\le1$, put
\[
 x_\tau=(y,\tau\xi),\quad
 J_s(\tau)=j_{i,s}(x_\tau,Q_{\tau\ell}(x_\tau)),
 \quad p_s=J_s(0).
\]
Equation~\eqref{kp:eq:model-jacobian} gives
$\partial_s^2p_s=0$ and $|p_s|+|\partial_sp_s|\le C$.

Let $g=(g_{ab})$ be the coordinate metric matrix and
define the weighted volume density of the coordinate
$y$-planes by
\[
 m(y,\xi)
 =w_i(y,\xi)\sqrt{\det(g_{ab}(y,\xi))_{a,b=1}^k},
 \quad m(y,0)=w_i(p).
\]
Write $\alpha=\alpha_s(|y|)$, and let $\nabla_\xi m$
denote the ordinary gradient in the $\xi$ variables.
Set
\[
 b_s=
 \frac{\alpha\nabla_\xi m(\alpha y,0)
       -\nabla_\xi m(y,0)}{w_i(p)}.
\]
Since $|x_\tau|^2=|y|^2+\tau^2|\xi|^2$, the tangential
block of $Df_s(x_\tau)$ has zero $\tau$-derivative at $0$.
The normal components and the plane variation give no
first-order contribution, by orthogonality and
\eqref{kp:eq:plane-first-derivative}.
Only the density $m$ contributes, and the quotient rule gives
\[
 \begin{aligned}
 \partial_\tau J_s(0)
 &=p_s\left.\frac{d}{d\tau}
   \frac{m(\alpha y,\alpha\tau\xi)}
        {m(y,\tau\xi)}\right|_{\tau=0}\\
 &=p_s\langle b_s,\xi\rangle.
 \end{aligned}
\]
Dots below denote $s$-derivatives. Direct differentiation gives
$\dot\alpha=\chi_r(|y|)\alpha^{1-k}$ and
$\ddot\alpha=(1-k)\chi_r(|y|)^2\alpha^{1-2k}$;
since $\alpha\ge1$, we have $|\dot\alpha|\le1$ and
$|\ddot\alpha|\le k-1$.
The derivatives of $m$ through order three are uniformly
bounded and $w_i(p)$ is bounded below.
Hence the product and chain rules give
$|b_s|+|\dot b_s|+|\ddot b_s|\le C(1+|y|+|y|^2)\le C$.
Since $\partial_s^2p_s=0$,
\[
 |\partial_s^2\partial_\tau J_s(0)|
 =\left|\left\langle
 2(\partial_sp_s)\dot b_s+p_s\ddot b_s,\xi
 \right\rangle\right|
 \le C|\xi|.
\]

For the second-order remainder, rescale coordinates by
$x=r\widehat x$. The map
\[
 \Phi_s(\widehat x)=r^{-1}f_s(r\widehat x)
 =(1+ks\chi(8|\widehat x|))^{1/k}\widehat x
\]
is independent of $r$, while
\[
D^q[g(r\,\cdot)]=r^q(D^qg)(r\,\cdot),
\quad
D^q[w_i(r\,\cdot)]=r^q(D^qw_i)(r\,\cdot),
\quad 0\le q\le4.
\]
Thus the entries of the source and target Gram matrices have uniformly
bounded derivatives through order four in the rescaled coordinates, and
their determinants and the weight denominator are bounded away from
zero. For a normal coordinate vector $\eta$ and a plane matrix $L$, put
$\widehat J_s(\eta,L)=j_{i,s}\bigl((y,r\eta),Q_L((y,r\eta))\bigr)$.
Differentiating the weighted Jacobian gives
$\bigl\|D_{(\eta,L)}^2\partial_s^2\widehat J_s\bigr\|\le C$, where
$D_{(\eta,L)}^2$ is the Hessian in $\eta$ and the entries of $L$.
Since $J_s(\tau)=\widehat J_s(\tau\xi/r,\tau\ell)$, this path has
constant velocity $v=(\xi/r,\ell)$, and consequently
\[
\left|\partial_\tau^2\partial_s^2J_s(\tau)\right|
=\left|D_{(\eta,L)}^2\partial_s^2\widehat J_s
(\tau\xi/r,\tau\ell)[v,v]\right|
\le C\bigl(r^{-2}|\xi|^2+|\ell|^2\bigr)
\le C\bigl(|\xi|/r+|\ell|\bigr)^2.
\]
Taylor's formula, $J_s(1)=j_{i,s}(x,Q)$, and
$\partial_s^2J_s(0)=0$ now imply
\begin{equation}\label{kp:eq:tangential-derivative-error}
 \begin{aligned}
 |\partial_s^2j_{i,s}(x,Q)|
 &\le|\partial_s^2\partial_\tau J_s(0)|
 +\tfrac12\sup_{0\le\tau\le1}
 |\partial_\tau^2\partial_s^2J_s(\tau)|\\
 &\le C\bigl(|\xi|+(|\xi|/r+|\ell|)^2\bigr)\\
 &\le C\bigl(u+r^{-2}u^2+a^2\bigr).
 \end{aligned}
\end{equation}
When $a(x,Q)$ is bounded away from zero, the same estimate
follows from the uniform bound
$|\partial_s^2j_{i,s}(x,Q)|\le C$.
By \eqref{kp:eq:finite-contact-errors} and Cauchy--Schwarz,
\[
 \int_{B_r}\bigl(u+r^{-2}u^2+a(x,T_xV)^2\bigr)\,d\|V\|
 =o(r^{k+2}).
\]
Apply \eqref{kp:eq:cutoff-convergence} with the test
$w_i\zeta_r(u+r^{-2}u^2+a^2)$.
For fixed $r$ and sufficiently large $i$, integrating
\eqref{kp:eq:tangential-derivative-error} twice in $s$
and using $\partial_sA_i(0,0)=0$ gives
\begin{equation}\label{kp:eq:expansion-bound}
 A_i(s,0)-\Mass_{w_i}(R_i)
 \le\omega(r)s^2r^{k+2},
 \quad \omega(r)\longrightarrow0.
\end{equation}
The same nonnegative error integral bounds $|\partial_s^2A_i(s,0)|$
for all $0\le s\le s_0$. The estimate is therefore uniform in $s$, and
$i$ can be chosen before $s$, without passing to a new subsequence for
each comparison map.

\noindent\textbf{Step 4: Inward motion.}\enspace
We view the uncut inward motion on the dilated current through
$f_s\circ G_t\circ f_s^{-1}$, since
\[
 f_s\circ G_t=(f_s\circ G_t\circ f_s^{-1})\circ f_s.
\]
Its infinitesimal field is $Y_s=(f_s)_*Y$. We compute the first-order
effect and bound the second-order remainder to prove
\eqref{kp:eq:inward-body}. On $\spt Y$, we have
$f_s(x)=\lambda x$ and
$Y_s(\lambda x)=\lambda r^2\rho_r(x)\partial_d$.
Let $P^\lambda$ be the coordinate matrix of the orthogonal
projection onto $Df_s(Q)$ at $\lambda x$.
First variation and the product rule give
\[
\partial_tA_i(s,0)
=r^2\int w_i j_{i,s}
\left[
 \lambda\rho_r
 \operatorname{div}_{Df_s(Q),w_i}\partial_d(\lambda x)
 +d\rho_r(P^\lambda\partial_d)
\right]d|R_i|(x,Q).
\]
Here $d\rho_r$ is evaluated at $x$.

Since $|\lambda x-x|\le Crs$, the bounded derivatives
of the metric and weight give
\[
|j_{i,s}-\lambda^k|\le Crs,\quad
\left|
 \operatorname{div}_{Df_s(Q),w_i}\partial_d(\lambda x)
 -\operatorname{div}_{Q,w_i}\partial_d(x)
\right|\le Crs.
\]
For the projection term, let the columns of $E$ be a
Euclidean orthonormal basis of the coordinate plane $Q$,
and let $e_d$ be the last coordinate vector.
Since $g e_d=e_d$,
\[
P^\lambda e_d
=E(E^Tg(\lambda x)E)^{-1}E^Te_d,
\quad
|E^Te_d|\le Ca(x,Q).
\]
The inverse matrices are uniformly bounded and change
by $O(rs)$. Together with $|\lambda-1|\le s$, this yields
$|j_{i,s}P^\lambda e_d-\lambda^{k+1}P^1e_d|\le Cs\,a(x,Q)$.
Subtract $\lambda^{k+1}\partial_tA_i(0,0)$ from the
first formula. The terms containing $\rho_r$ are bounded
by $Csr^3$; those containing $d\rho_r$ are bounded by
$Csra$, because $|d\rho_r|\le C/r$. Thus
\[
\begin{aligned}
|\partial_tA_i(s,0)+\lambda^{k+1}B_i|
&\le Cs\left[
 r^3\|R_i\|(B_{r/2})
 +r\int_{B_{r/2}}a(x,T_xR_i)\,d\|R_i\|
\right]\\
&\le s\omega(r)r^{k+2},
\quad \omega(r)\longrightarrow0.
\end{aligned}
\]
We used $\partial_tA_i(0,0)=-B_i$ from
\eqref{kp:eq:boundary-force}. The last inequality follows
from \eqref{kp:eq:finite-contact-errors},
Cauchy--Schwarz, and \eqref{kp:eq:cutoff-convergence},
for fixed small $r$ and sufficiently large $i$.
More explicitly, the corresponding limit error is bounded by
\[
 Cr^{k+3}+Cr\,r^{k/2}
 \left(\int_{B_r}a(x,T_xV)^2\,d\|V\|\right)^{1/2}
 =o(r^{k+2}).
\]
The cutoff $\zeta_r$ transfers this bound to $R_i$ without requiring
convergence on the boundary of a metric ball.

For the quadratic remainder, on $\spt Y$ we have
$F_t:=f_s\circ G_t=\lambda(x+tY)$,
$DF_t=\lambda(\id+tDY)$, $|Y|\le Cr^2$, and $|DY|\le Cr$.
For an orthonormal basis $e_1,\ldots,e_k$ of $Q$, set
\[
H_{ab}(t)=
\langle DF_t(x)e_a,DF_t(x)e_b\rangle_{F_t(x)},
\quad
W(t)=\frac{w_i(F_t(x))}{w_i(x)}.
\]
Primes denote $t$-derivatives. Since $F_t$ and $DF_t$ are affine in
$t$, the product and chain rules give
\[
\begin{aligned}
|W'|&\le C|Y|,&
|W''|&\le C|Y|^2,\\
|H'|&\le C(|Y|+|DY|),&
|H''|&\le C(|Y|^2+|Y||DY|+|DY|^2).
\end{aligned}
\]
The matrices $H,H^{-1}$ are uniformly bounded because
$|tDY|<1/2$. Differentiating $\mathcal J=W\sqrt{\det H}$ twice
therefore gives
$|\mathcal J''|\le C\bigl(|W''|+|W'||H'|+|H''|+|H'|^2\bigr)\le Cr^2$.
Integrating over $\spt Y\subset B_{r/2}$ yields
$|\partial_t^2A_i(s,t)|\le Cr^2\|R_i\|(B_{r/2})\le Cr^{k+2}$.
Taylor's formula now gives
\begin{equation}\label{kp:eq:inward-body}
A_i(s,t)-A_i(s,0)
\le -t\lambda^{k+1}B_i
     +st\omega(r)r^{k+2}
     +Ct^2r^{k+2}.
\end{equation}

\noindent\textbf{Step 5: The boundary layer.}\enspace
We estimate the cost of cutting off the inward motion in the collar $\{0<d<\varepsilon\}$; the result is \eqref{kp:eq:inward-layer}.
Fix $i,s,t$. Take $\varepsilon\downarrow0$ through values
with $\|R_i\|(\{d=\varepsilon\})=0$.
In $0<d<\varepsilon$, put $\tau=td/\varepsilon$ and
$A=D(f_s\circ G_\tau)_x$. The chain rule gives
\[
 D(f_s\circ F_{t,\varepsilon})v
 =Av+\frac{\lambda tr^2\rho_r(x)}{\varepsilon}
       \langle\nabla d,v\rangle\partial_d,
\]
where the final $\partial_d$ is evaluated at $f_sG_\tau x$.

If $\pi_Q\nu\ne0$, choose an orthonormal basis
$e_1,\ldots,e_k$ of $Q$ with
$e_1=\pi_Q\nu/|\pi_Q\nu|$. Since $\nabla d=-\nu$, the
additional rank-one term acts only on $e_1$. With
$c=\lambda tr^2\rho_r(x)/\varepsilon$, the exterior-product
formula and $|\partial_d|=1$ give
\[
\begin{aligned}
 J_k\bigl(D(f_s\circ F_{t,\varepsilon})|Q\bigr)
 &=\bigl|
   (Ae_1-c|\pi_Q\nu|\partial_d)
   \wedge Ae_2\wedge\cdots\wedge Ae_k
   \bigr|\\
 &\le J_k(A|Q)
   +\frac{\lambda tr^2\rho_r(x)}{\varepsilon}
    |\pi_Q\nu|\,\|A\|^{k-1},
\end{aligned}
\]
where $\|A\|$ is the operator norm. If $\pi_Q\nu=0$, the
additional term vanishes, so the estimate remains valid.

On $\{\rho_r>0\}$, Step~2 gives $G_\tau(x)\in\spt\rho_r$, where
$Df_s=\lambda\id$. Hence $A=\lambda(\id+\tau DY)$ there, and
\[
 \frac{w_i(f_sG_\tau x)}{w_i(x)}
 \lambda\|A\|^{k-1}\le\lambda^k+Cr(s+t).
\]
Indeed, the metric and weight change by $O(rs)$ from $x$
to $\lambda x$, while $|\tau DY|\le Crt$ and
$|\tau Y|\le Cr^2t$. On $\{\rho_r=0\}$ the additional Jacobian term
vanishes.

Put $E_\varepsilon=\{0<d<\varepsilon\}$.
Because $R_i$ is canonical, $\|R_i\|(E_\varepsilon)\to0$.
The maps $f_s\circ G_\tau$ and their differentials, with
$\tau$ held fixed, are uniformly bounded for $0\le\tau\le t$;
thus the integral of their weighted Jacobian over $E_\varepsilon$
tends to zero. On $\{d>\varepsilon\}$, $F_{t,\varepsilon}=G_t$, so
\[
 \int_{\{d>\varepsilon\}}
 w_i(f_s(G_t(x)))
 J_k\bigl(D(f_s\circ G_t)_x|T_xR_i\bigr)\,
 d\|R_i\|(x)
 \longrightarrow A_i(s,t).
\]
For the additional Jacobian term, use
$|\pi_{T_xR_i}\nu|\le1$ and the preceding coefficient
estimate to obtain
\[
\begin{aligned}
 &\frac{\lambda tr^2}{\varepsilon}
 \int_{E_\varepsilon}
 w_i(f_s(G_\tau(x)))\rho_r(x)
 |\pi_{T_xR_i}\nu|\,\|A\|^{k-1}\,d\|R_i\|(x)\\
 &\quad\le
 \frac{tr^2}{\varepsilon}
 \bigl[\lambda^k+Cr(s+t)\bigr]
 \int_{E_\varepsilon}\rho_r w_i\,d\|R_i\|.
\end{aligned}
\]
Apply the continuous-extension version of
\eqref{kp:eq:weighted-collar-localization}
in Lemma~\ref{kp:lem:boundary-force},
with $\widetilde b=\rho_r$ and
$b=\rho_r|_{\partial M}$. This gives
\[
 \limsup_{\varepsilon\downarrow0}
 \frac1\varepsilon
 \int_{E_\varepsilon}\rho_r w_i\,d\|R_i\|
 \le\int_{\partial M}\rho_r w_i\,d\sigma_i
 =\frac{B_i}{r^2}.
\]
Combining these estimates and using $B_i\le Cr^{k+2}$,
we obtain
\begin{equation}\label{kp:eq:inward-layer}
\begin{aligned}
 \limsup_{\varepsilon\downarrow0}
 \Mass_{w_i}(S_i(s,t,\varepsilon))
 &\le A_i(s,t)+t\bigl[\lambda^k+Cr(s+t)\bigr]B_i\\
 &=A_i(s,t)+t\lambda^kB_i+Crt(s+t)B_i\\
 &\le A_i(s,t)+t\lambda^kB_i+Ct(s+t)r^{k+3}.
\end{aligned}
\end{equation}

\noindent\textbf{Step 6: The admissible endpoint.}\enspace
We add the three estimates, choose the parameters, and verify that the resulting current belongs to $\mathcal C_i$.
Since $\lambda-1\ge s/2$ and $B_i\ge cr^{k+2}$,
$t\lambda^kB_i-t\lambda^{k+1}B_i=-t\lambda^k(\lambda-1)B_i\le-cst\,r^{k+2}$,
after decreasing $c$.
Set $t=\eta s$, where $0<\eta\le1$.
Adding \eqref{kp:eq:expansion-bound},
\eqref{kp:eq:inward-body}, and \eqref{kp:eq:inward-layer},
and enlarging $\omega(r)$ to absorb $Cr$, gives
\begin{equation}\label{kp:eq:finite-area-bounds}
 \limsup_{\varepsilon\downarrow0}
 [\Mass_{w_i}(S_i(s,\eta s,\varepsilon))-\Mass_{w_i}(R_i)]
 \le s^2r^{k+2}[(1+\eta)\omega(r)-c\eta+C\eta^2].
\end{equation}

It remains to control the additional square term in $F_i$
and to verify that the new relative cycle belongs to
$\mathcal C_i$.

Write $\Psi=f_s\circ F_{t,\varepsilon}$,
$S=S_i(s,t,\varepsilon)$, and
$J(x)=J_k(D\Psi_x|T_xR_i)$.
The reverse triangle inequality applied to the preceding
exterior products gives, after integration,
\[
\begin{aligned}
 \int|J-1|\,d\|R_i\|
 &\le C(s+rt)\|R_i\|(B_{r/2})
   +\frac{Ctr^2}{\varepsilon}
    \int_{\{0<d<\varepsilon\}}\rho_r\,d\|R_i\|\\
 &\le C_r(s+t).
\end{aligned}
\]
The last inequality follows from
\eqref{kp:eq:weighted-collar-localization},
$w_i\ge1/2$, and $B_i\le Cr^{k+2}$, for sufficiently
small $\varepsilon$ at fixed $i$.
In fact, the rescaled weighted collar integral is at most
$2B_i/r^2$ for sufficiently small $\varepsilon$, because
$B_i>0$. Thus $C_r$ is independent of $i$; only the required
upper bound on $\varepsilon$ may depend on the already fixed $i$.
Moreover,
$\sup_x\dist(\Psi(x),x)\le C(rs+r^2t)$.
Since $\Psi$ is injective, the area formula therefore gives
\begin{equation}\label{kp:eq:inward-moments}
\begin{aligned}
 \left|\int f\,d(\|S\|-\|R_i\|)\right|
 &=\left|\int[f(\Psi(x))J(x)-f(x)]\,d\|R_i\|(x)\right|\\
 &\le \|f\|_\infty\int|J-1|\,d\|R_i\|
   +\|\nabla f\|_\infty
      \int\dist(\Psi(x),x)\,d\|R_i\|(x)\\
 &\le C_r(s+t)\|f\|_{C^1}.
\end{aligned}
\end{equation}

We next estimate the relative flat distance. Although $R_i$ has finite mass, finiteness of $\Mass(\partial R_i)$ has not been established.
To apply the homotopy formula for integral currents, we therefore choose a representative with finite boundary mass. By \cite[Lemma 3.1]{GLWZ2021}, choose an integral representative $\widetilde R_i$ of $[R_i]$ with $\Mass(\widetilde R_i)\le\Mass(R_i)+1$.
For $0\le v\le1$, set
$H(v,x)=x+vtq_\varepsilon(d(x))Y(x)$.
Then
\[
 \partial_vH=tq_\varepsilon Y,\quad
 D_xH(e)=e+vtq_\varepsilon D_eY
              +vt\,dq_\varepsilon(e)Y.
\]
The last term vanishes when wedged with $\partial_vH$.
Consequently the factor $q_\varepsilon'$, although of size
$1/\varepsilon$, does not enter the swept-mass bound.
Since $|Y|\le Cr^2$ and $|DY|\le Cr$, the current
swept out by $H$ has mass at most
$Ctr^2\Mass(\widetilde R_i)$.
Push this homotopy forward by $f_s$, whose Lipschitz
constant is uniformly bounded, and add the homotopy
$u\mapsto f_u$, $0\le u\le s$, applied to
$\widetilde R_i$.
The latter sweeps mass at most
$Crs\Mass(\widetilde R_i)$.
Both homotopies keep the boundary of the current in
$\partial M$. The relative homotopy formula thus gives
\begin{equation}\label{kp:eq:inward-flat}
 \Flat([S]-[R_i])
 \le C(tr^2+sr)\Mass(\widetilde R_i)
 \le C_r(s+t).
\end{equation}

By \eqref{kp:eq:auxiliary-expansion},
\eqref{kp:eq:test-bound}, and
\eqref{kp:eq:inward-moments},
\[
 F_i(S)-F_i(R_i)
 \le \Mass_{w_i}(S)-\Mass_{w_i}(R_i)
       +\frac{C_r}{i}(s+t)^2.
\]
Choose $\eta>0$, then $r>0$ sufficiently small, and then $i$ sufficiently large so that
\[
 C\eta<\frac c4,\quad
 (1+\eta)\omega(r)<\frac{c\eta}{4},\quad
 \frac{C_r}{i}(1+\eta)^2
       <\frac{c\eta}{4}r^{k+2}.
\]
The index $i$ is also chosen large enough for the
preceding estimates.
The order is essential: after $r$ is fixed, the term $C_r/i$ can
be made negligible; after $i$ is fixed, its possibly very small
admissibility scale $\delta_i$ is accommodated by choosing $s$.
Set $t=\eta s$ and choose $0<s<s_0$ with
$C_r(1+\eta)s<\delta_i$.
Equation~\eqref{kp:eq:finite-area-bounds} now yields
\[
 \limsup_{\varepsilon\downarrow0}
 \bigl[F_i(S_i(s,\eta s,\varepsilon))-F_i(R_i)\bigr]
 \le-\frac{c\eta}{4}s^2r^{k+2}<0.
\]
Fix a sufficiently small $\varepsilon$ so that
$F_i(S)<F_i(R_i)$. Then
$\Mass(S)\le F_i(S)<F_i(R_i)\le\Mass(T_i)$.
Also, $\spt(S-R_i)\subset B_{r/2}\subset K$, and
\eqref{kp:eq:inward-flat} gives
$\Flat([S]-[R_i])<\delta_i$.
Thus $[S]$ can be appended to the sequence reaching
$[R_i]$ in \eqref{kp:eq:reachable}, so $[S]\in\mathcal C_i$.
Only this final endpoint is appended. The homotopies above estimate
its flat distance; they are not required to stay below the mass
threshold at every intermediate time.
This contradicts the minimality of $R_i$ in
Lemma~\ref{kp:lem:comparison-minimizers}.
\end{proof}

\begin{proof}[Proof of Proposition~\ref{kp:local-trace}]
Fix $\Omega'\Subset\Omega\Subset U$ and write the restriction of
$\sigma$ on $\Omega'$ as $q\,\|V\|\res(\Omega'\cap\partial M)+\sigma_s$.
Put $h(x)=\operatorname{tr}_{T_xV}\SFF_{\partial M}$.
Theorem~\ref{kp:thm:nonnegative-trace} gives $h=q\ge0$ at
$\|V\|$-almost every point of $\Omega'\cap\partial M$.
Lemma~\ref{kp:lem:positive-trace} excludes $h>0$ on a set of positive
$\|V\|$-measure there. Thus $h=q=0$ almost everywhere on
$\Omega'\cap\partial M$, and the restriction of $\sigma$ to $\Omega'$
is concentrated on a $\|V\|$-null set. Along a countable exhaustion of
$U$ this proves \eqref{kp:eq:local-trace}. The union of the corresponding
null sets is $\|V\|$-null and carries $\sigma$, proving singularity on $U$.
\end{proof}

\subsection{Geometric consequences}
\label{sec:consequences}

\begin{proof}[Proof of Theorem~\ref{thm:kflat}]
The varifold $V$ is integral by Corollary~\ref{loc:cor:integral}. Choose $\rho_p>0$ below the small-annulus radius at each $p$,
and a finite cover $M=\bigcup_{\ell=1}^N B_{\rho_\ell}(p_\ell)$,
where $\rho_\ell=\rho_{p_\ell}$. Each punctured ball is covered by
annuli in which $V$ is almost minimizing:
\[
 B_{\rho_\ell}(p_\ell)\setminus\{p_\ell\}
 =
 \bigcup_{j=1}^{\infty}
 \An(p_\ell;2^{-j}\rho_\ell,\rho_\ell).
\]
Put $E=\{p_1,\ldots,p_N\}$. Proposition~\ref{kp:local-trace} gives
the trace identity almost everywhere on $M\setminus E$. Since $V$ is
integral, $\|V\|(E)=0$, which proves the identity on $M$.

On $M\setminus E$, the measures in Theorem~\ref{thm:local}\textup{(ii)}
agree on overlaps by uniqueness. The countable annular cover and
\textup{(iii)} thus give $\delta V=\nu\,\sigma$ there, with
$\sigma\ge0$ and $\sigma\perp\|V\|$. To extend this formula across
$E$, choose a smooth cutoff $0\le\chi_r\le1$ equal to one on $B_r(E)$,
supported in $B_{2r}(E)$, with $|\nabla\chi_r|\le C/r$.
The monotonicity estimate in the proof of Corollary~\ref{loc:cor:integral}
gives $\|V\|(B_{2r}(E))=O(r^k)$, since $E$ is finite. Hence for every
smooth field $X$,
\[
 |\delta V(\chi_r X)|
 \le C_X(1+r^{-1})\|V\|(B_{2r}(E))
 =O_X(r^{k-1})\longrightarrow0.
\]
First take $X$ to be a smooth extension of $\nu$ from $\partial M$.
Then $\int(1-\chi_r)\,d\sigma=\delta V((1-\chi_r)X)\to\delta V(X)$,
so Fatou's lemma gives $\sigma(\partial M\setminus E)<\infty$.
Extend $\sigma$ by zero on $E$. For arbitrary $X$, dominated
convergence now gives
\[
 \delta V(X)=\lim_{r\downarrow0}\delta V((1-\chi_r)X)
 =\int_{\partial M}g(X,\nu)\,d\sigma .
\]
The extension remains singular with respect to $\|V\|$, and testing
with smooth multiples of the normal determines $\sigma$ uniquely.
\end{proof}

\begin{proof}[Proof of Corollary~\ref{cor:zeroboundary}]
By Theorem~\ref{thm:kflat}, almost every boundary-mass point has
$T_xV\subset T_x\partial M$ and zero trace, contradicting either
alternative in \eqref{eq:stricttrace}. Thus $\|V\|(\partial M)=0$.
\end{proof}

\begin{proof}[Proof of Corollary~\ref{cor:euclidean-interior}]
Using a smooth boundary collar, extend
$x-\langle x,\nu\rangle\nu$ on $\partial\Omega$ to a smooth vector
field $Y$ on $\overline\Omega$. For every $k$-plane
$P\subset T_x\partial\Omega$,
\[
 \diver_PY(x)
 =k-\langle x,\nu(x)\rangle
       \operatorname{tr}_P\SFF_{\partial\Omega}(x).
\]
By Theorem~\ref{thm:kflat}, $\diver_{T_xV}Y(x)=k$ for
$\|V\|$-almost every $x\in\partial\Omega$. Fix this extension and
put $C=\max_{(x,P)\in G_k(M)}(-\diver_PY(x))_+$. Since $Y$ is
tangent to $\partial\Omega$, free-boundary stationarity gives
\[
 0=\delta V(Y)
 =k\|V\|(\partial\Omega)+\int_\Omega\diver_{T_xV}Y\,d\|V\|
 \ge k\|V\|(\partial\Omega)-C\|V\|(\Omega).
\]
Adding $k\|V\|(\Omega)$ to both sides of
$k\|V\|(\partial\Omega)\le C\|V\|(\Omega)$ gives the claim with
$c_{\Omega,k}=k/(k+C)$.
\end{proof}

The barrier argument uses the sign convention fixed in \eqref{eq:boundaryform}.

\begin{proposition}
\label{prop:barrier}
Suppose $H_{\partial M}>0$. Then every smooth almost properly embedded free-boundary minimal hypersurface $\Sigma\subset M$ is properly embedded, that is, $\Sigma\cap\partial M=\partial\Sigma$. Here almost properly embedded has the meaning of \cite[Definition 2.6]{LZ2021}: the hypersurface is smoothly embedded, $\partial\Sigma\subset\partial M$, and interior contact with $\partial M$ is allowed.
\end{proposition}

\begin{proof}
Let $r$ be the inward distance to $\partial M$. Since $\nabla r=-\nu$
along the boundary,
\begin{equation}\label{eq:distancehessian}
 \nabla_M^2r(v,w)=-\SFF_{\partial M}(v,w)
 \quad(v,w\in T\partial M).
\end{equation}
Suppose that an interior point $p$ of $\Sigma$ lies on $\partial M$.
The function $r|_\Sigma\ge0$ attains its minimum at $p$, so
$T_p\Sigma=T_p\partial M$. Since $\Sigma$ is a smooth minimal
hypersurface near $p$,
\[
 0\le\Delta_\Sigma(r|_\Sigma)(p)
 =\operatorname{tr}_{T_p\Sigma}\nabla_M^2r
 =-H_{\partial M}(p)<0,
\]
a contradiction. This is the smooth interior-contact case of the
maximum principle; compare \cite[Theorem 1]{White2010}.
\end{proof}

\begin{proof}[Proof of Corollary~\ref{cor:hypersurface}]
For $3\le n\le7$, the mod-two free-boundary regularity theorem \cite[Theorem 4.6]{GLWZ2021}, which is the mod-two form of the regularity theory of Li--Zhou \cite{LZ2021}, gives the representation \eqref{eq:hypersurfacerepresentation}, in which $N$ and the $a_i$ are positive integers and each $\Sigma_i$ is a smooth, compact, connected, almost properly embedded free-boundary minimal hypersurface. Proposition~\ref{prop:barrier} makes every $\Sigma_i$ properly embedded. Its boundary is a smooth $(n-2)$-dimensional submanifold, possibly empty. Consequently,
\[
 \|V\|(\partial M)
 =\sum_{i=1}^{N}a_i\HH^{n-1}(\Sigma_i\cap\partial M)
 =\sum_{i=1}^{N}a_i\HH^{n-1}(\partial\Sigma_i)=0.
\]
\end{proof}

We conclude with the proof of Corollary~\ref{cor:dirichlet}, which
applies Corollary~\ref{cor:euclidean-interior} to a varifold limit.

\begin{proof}[Proof of Corollary~\ref{cor:dirichlet}]
Fix $\Lambda$, $N$, and $R_*$, and let $\mathcal S$ be the class of
submanifolds $\Sigma$ described in the statement. We may assume that
$\mathcal S$ is nonempty.

Let $\Sigma\in\mathcal S$. In \eqref{eq:smoothvariation},
$\mathbf H_\Sigma=0$ and the conormal $\eta_\Sigma$ is normal to
$\partial\Omega$ along $\partial\Sigma$, so $|\Sigma|$ is stationary
with free boundary. Minimality also gives $\Delta_\Sigma|x|^2=2k$ on
$\Sigma$. If $\partial\Sigma$ were empty, integrating this identity
over $\Sigma$ would give $0=2k\,\HH^k(\Sigma)$, which is impossible.
Hence $\partial\Sigma\ne\emptyset$. The density of $|\Sigma|$ at a
point of $\partial\Sigma\subset\partial\Omega$ is $1/2$, so the
monotonicity formula at boundary points \cite[Theorem~2.3]{LZ2021}
gives a constant $c_0>0$, depending only on $\Omega$ and $k$, such that
\begin{equation}\label{eq:monotone-bounds}
 \HH^k(\Sigma)\ge c_0
 \quad\text{for every }\Sigma\in\mathcal S.
\end{equation}

Set $d_*=\inf_{\Sigma\in\mathcal S}d_\Sigma$ and suppose that $d_*=0$.
Choose $\Sigma_i\in\mathcal S$ with $d_{\Sigma_i}\to0$. Since
$\HH^k(\Sigma_i)\le\Lambda$, a subsequence of the varifolds
$|\Sigma_i|$, not relabeled, converges to a $k$-varifold $V$ in
$\overline\Omega$. The first variation on a fixed smooth vector field is
continuous under varifold convergence, so $V$ is stationary with free
boundary. By \eqref{eq:monotone-bounds},
$\|V\|(\overline\Omega)=\lim_i\HH^k(\Sigma_i)\ge c_0$, so $V\ne0$.

We show that $V$ has property $\mathcal P_N$ at scale $R_*$; compare
\cite[Lemma~3.9]{PartII}. Almost minimization in a fixed relatively
open set $U\subset\overline\Omega$ is preserved under varifold
convergence. Indeed, let $W_i\to W$, where each $W_i$ is almost
minimizing modulo two in $U$, and let $\eps>0$. Choose $i$ with
$\dV(W_i,W)<\eps/2$, and then $\delta>0$ and
$T\in\AMset(U;\eps/2,\delta;\Flat)$ with $\dV(|T|,W_i)<\eps/2$. Then
$T\in\AMset(U;\eps,\delta;\Flat)$ and $\dV(|T|,W)<\eps$. Now fix
$p\in\overline\Omega$ and radii as in \eqref{eq:PN}. Each $|\Sigma_i|$
is almost minimizing in $\An(p;s_j,r_j)$ for some $j=j(i)$, and one
value of $j$ occurs for infinitely many $i$. Hence $V$ is almost
minimizing in the corresponding annulus.

Property $\mathcal P_N$ implies that $V$ is almost minimizing in small
relative annuli. Otherwise there would be a point
$p\in\overline\Omega$ such that, for every $\rho>0$, $V$ is not almost
minimizing in some $\An(p;s,r)$ with $0<s<r<\rho$. Choosing such an
annulus with $r_1<R_*$, and then successively with $r_{j+1}<s_j/12$,
would give $N$ annuli as in \eqref{eq:PN} in none of which $V$ is
almost minimizing.

Finally, $\|V\|(\Omega)=0$. Let $\varphi\in C_c(\Omega)$. Its support
has positive distance from $\partial\Omega$, whereas every point of
$\Sigma_i$ has distance at most $d_{\Sigma_i}$ from $\partial\Omega$.
Hence $\varphi$ vanishes on $\Sigma_i$ for all large $i$, and
\[
 \int\varphi\,d\|V\|
 =\lim_{i\to\infty}\int_{\Sigma_i}\varphi\,d\HH^k=0.
\]
Thus $V$ is a nonzero $k$-varifold that is stationary with free
boundary and almost minimizing modulo two in small relative annuli,
with $\|V\|(\Omega)=0$. This contradicts
Corollary~\ref{cor:euclidean-interior}. Therefore $d_*>0$, which proves
\eqref{eq:depth}.

It remains to prove \eqref{eq:dirichlet}. Fix $\Sigma\in\mathcal S$,
write $d=d_\Sigma$, and choose $q\in\Sigma$ with
$\dist(q,\partial\Omega)=d$. Put
\[
 u(x)=\bigl(d-|x-q|\bigr)_+ .
\]
Since $\partial\Sigma\subset\partial\Omega$, the Lipschitz function
$u$ vanishes on $\partial\Sigma$, and $u(q)=d>0$. Its restriction to
$\Sigma$ is therefore an admissible test function for
$\lambda_1^D(\Sigma)$. The tangential gradient satisfies
$|\nabla_\Sigma u|\le1$ almost everywhere, whence
\[
 \int_\Sigma|\nabla_\Sigma u|^2\,d\HH^k\le\HH^k(\Sigma).
\]
For $0<t<d$ the Euclidean ball $B_t(q)$ is contained in $\Omega$ and
hence does not meet $\partial\Sigma$, so the monotonicity formula
\cite[Section~5]{Allard1972} at the smooth point $q$ gives
$\HH^k(\Sigma\cap B_t(q))\ge\omega_kt^k$. Consequently,
\begin{align*}
 \int_\Sigma u^2\,d\HH^k
 &=2\int_0^{d}(d-t)\,\HH^k(\Sigma\cap B_t(q))\,dt\\
 &\ge2\omega_k\int_0^{d}(d-t)\,t^k\,dt
  =\frac{2\omega_k}{(k+1)(k+2)}\,d^{\,k+2}.
\end{align*}
The quotient of these two integrals bounds $\lambda_1^D(\Sigma)$ from
above, which gives the first inequality in \eqref{eq:dirichlet}. The
second follows from $\HH^k(\Sigma)\le\Lambda$ and $d_\Sigma\ge d_*$.
\end{proof}

The hypothesis $\mathcal P_N$ in Corollary~\ref{cor:dirichlet} cannot
be replaced by the area bound alone, even in the unit ball
$B\subset\R^3$. Let $\Sigma_i\subset\overline B$ be the smooth embedded
genus-zero free-boundary minimal surfaces of
\cite[Corollary~1.5]{KS2024}, for which $|\Sigma_i|\to|\partial B|$
as varifolds. Their areas converge to $4\pi$, and the surfaces
converge to $\partial B$ in Hausdorff distance
\cite[Corollary~2.27(1)]{KS2024}. Hence
\[
 m_i:=\max_{x\in\Sigma_i}(1-|x|^2)
 \le 2\max_{x\in\Sigma_i}(1-|x|)
 \le 2d_H(\Sigma_i,\partial B)\longrightarrow0.
\]
On $\Sigma_i$, the function
$u_i=1-|x|^2$ satisfies $-\Delta_{\Sigma_i}u_i=4$ by minimality, and
$u_i=0$ on $\partial\Sigma_i\subset\partial B$. Let $v_i>0$ be a first
Dirichlet eigenfunction of $\Sigma_i$. Green's identity gives
\[
 4\int_{\Sigma_i}v_i\,d\HH^2
 =\lambda_1^D(\Sigma_i)\int_{\Sigma_i}u_iv_i\,d\HH^2
 \le\lambda_1^D(\Sigma_i)\,m_i\int_{\Sigma_i}v_i\,d\HH^2.
\]
Hence $d_{\Sigma_i}\le m_i\to0$ and
$\lambda_1^D(\Sigma_i)\ge4/m_i\to\infty$, although the areas are
bounded. By Corollary~\ref{cor:dirichlet}, for no fixed $N$ and $R_*$
do all the varifolds $|\Sigma_i|$ have property $\mathcal P_N$ at
scale $R_*$.

\subsection{Examples}
\label{sec:examples}

Two examples show that the hypotheses of Theorem~\ref{thm:kflat} and Corollary~\ref{cor:zeroboundary} are needed: the relative comparison class in the first, and the strictness of the trace condition in the second.

\begin{remark}
\label{rem:relative-almost-minimizing}
Definition~\ref{def:almostminimizing} compares relative cycles using
canonical mass and relative flat distance. Its limiting varifold has
ordinary mass, including any boundary mass. Absolute almost minimization
instead compares currents with zero absolute boundary, without quotienting
out boundary-supported currents. The next example shows why this distinction
is essential in Theorem~\ref{thm:kflat}.
\end{remark}

\begin{example}
\label{ex:spherical-boundary}
Let
\[
 M=\{x\in S^3\subset\R^4:x_4\ge-1/2\},
 \quad \Gamma=\partial M,
\]
where $S^3$ is the unit round sphere, and put $W=|\Gamma|$, an integral $2$-varifold. Write $x=(\sin\theta\,\omega,\cos\theta)$ with $\omega\in S^2$. The induced metric is
\[
 g=d\theta^2+\sin^2\theta\,g_{S^2},
 \quad \Gamma=\{\theta=\theta_0\},\quad \theta_0=2\pi/3.
\]
For the outward normal $\nu=\partial_\theta$,
$\SFF_\Gamma=\cot\theta_0\,g_\Gamma=-g_\Gamma/\sqrt3$. Hence
$\operatorname{tr}_{T\Gamma}\SFF_\Gamma=-2/\sqrt3$ and
$\|W\|(\partial M)=4\pi\sin^2\theta_0=3\pi$.
Nevertheless $\delta W(X)=\int_\Gamma\diver_\Gamma X\,d\HH^2=0$
for every vector field tangent to $\Gamma$, by the divergence theorem.

Let $T=\llbracket\Gamma\rrbracket$ modulo two, so $\partial T=0$ and $W=|T|$. We prove the stronger local minimizing property: for every $p\in M$
there is $r_p>0$ such that every absolute mod-two cycle $S$ satisfies
\begin{equation}\label{ex:eq:absolute-local-minimum}
 \spt(S-T)\subset U_p:=M\cap B_{r_p}(p)
 \quad\Longrightarrow\quad\Mass(S)\ge\Mass(T).
\end{equation}

Balls are taken in a smooth extension, as in Section~\ref{sec:notation}.

For $p\in\Int M$, take $U_p$ disjoint from $\Gamma$;
$S=T$ near $\Gamma$ proves \eqref{ex:eq:absolute-local-minimum}.
For $p\in\Gamma$, use the collar
\[
 \mathcal C=\{\pi/2\le\theta\le\theta_0\}\subset M,\quad
 \pi:\mathcal C\longrightarrow\Gamma,\quad
 \pi(\theta,\omega)=(\theta_0,\omega).
\]
Since $\sin\theta\ge\sin\theta_0$ on $\mathcal C$, the radial
projection $\pi$ is $1$-Lipschitz: it kills radial vectors and multiplies
the lengths of tangential vectors by $\sin\theta_0/\sin\theta\le1$.

Choose $U_p\subset\mathcal C$ with $\pi(U_p)$ contained in a
coordinate disk whose closure is proper in $\Gamma$. For $S$ as in
\eqref{ex:eq:absolute-local-minimum}, $\spt S\subset\Gamma\cup U_p$.
Then $\pi_\#(S-T)$ is a top-dimensional cycle on the connected
surface $\Gamma$. Constancy makes it either zero or the whole
$\Gamma$ modulo two; support in a proper disk excludes the latter
\cite[4.1.31(2) and 4.2.26]{Federer1969}.
Since $\pi_\#T=T$, the area formula gives
\[
 \Mass(S)\ge\Mass(\pi_\#S)=\Mass(T)=3\pi,
\]
proving \eqref{ex:eq:absolute-local-minimum} also at boundary points.

By \eqref{ex:eq:absolute-local-minimum}, every absolute comparison
sequence in $\An(p;s,r)$, $0<s<r<r_p$, ends with mass at least
$\Mass(T)$, regardless of flat step size or mass allowance.
Thus the constant sequence $T$ verifies absolute almost minimization
despite negative boundary trace and positive boundary mass.
\end{example}

The next example satisfies relative almost minimization, and its
ambient manifold satisfies the homological hypothesis
$H_1(M,\partial M;\Zt)=0$ of the existence theorem in \cite{PartII},
yet the varifold concentrates on a totally geodesic boundary. A
non-strict trace condition therefore does not suffice.

\begin{proposition}
\label{prop:boundary-concentration}
There is a smooth Riemannian metric on the closed three-ball $M$ such that $\partial M$ is totally geodesic and the integral $2$-varifold $V=|\partial M|$ is stationary with free boundary and almost minimizing modulo two in $M$. In particular, it is almost minimizing in small relative annuli, but
\begin{equation}\label{ex:eq:boundary-concentration}
 \|V\|(\partial M)=\|V\|(M)=4\pi.
\end{equation}
Moreover, $H_1(M,\partial M;\Zt)=0$.
\end{proposition}

\begin{proof}
\noindent\textbf{The metric and its minimal slices.}

Choose the smooth step function
\[
 \theta(s)=\frac{e^{-1/s}}{e^{-1/s}+e^{-1/(1-s)}}\quad(0<s<1),
 \quad \theta|_{(-\infty,0]}=0,\quad\theta|_{[1,\infty)}=1.
\]

For $0<t\le2$, set
\begin{equation}\label{ex:eq:ball-warping}
 f(t)=1+\bigl(1-\theta(t)\bigr)e^{-1/t^2}\sin^2(\pi/t)
          +(1-t)\theta(t-1),
 \quad f(0)=1.
\end{equation}
Collapse $\{2\}\times S^2$ in $[0,2]\times S^2$ to a point and equip
the resulting ball $M$ with $g=dt^2+f(t)^2g_{S^2}$,
where $g_{S^2}$ is unit round. Here $f\ge1$ on $[0,1]$,
$f(t)\ge2-t>0$ on $(1,2)$, and $f(2)=0$.
The identity $f(2-r)=r+(1-r)\theta(r)$ makes $f(2-r)-r$ flat
at $r=0$, ensuring smoothness at the pole. Likewise $f-1$ is
flat at $t=0$, so
\[
 \SFF_{\partial M}=-\frac{f'(0)}{f(0)}g|_{\partial M}=0,
 \quad\HH^2(\partial M)=4\pi.
\]
Since $M\cong B^3$, the homology assertion follows.

For $i\ge2$, put $t_i=1/i$, $\Sigma_i=\{t_i\}\times S^2$, and let $T_i$ be its relative mod-two current. Equation~\eqref{ex:eq:ball-warping} gives
\begin{equation}\label{ex:eq:slice-derivatives}
 f(t_i)=1,\quad f'(t_i)=0,\quad
 f''(t_i)=2\pi^2i^4e^{-i^2}\bigl(1-\theta(1/i)\bigr)>0.
\end{equation}
Thus each $\Sigma_i$ is minimal of area $4\pi$, and its second
variation is
\[
 \delta^2\operatorname{Area}_{\Sigma_i}(u,u)
 =\int_{S^2}\bigl(|\nabla u|^2+2f''(t_i)u^2\bigr)\,dA>0
 \quad(u\ne0).
\]

\noindent\textbf{Strict local minimality in relative flat distance.}

For each $i$, we claim that some $\rho_i>0$ makes every relative
mod-two cycle $T$ satisfy
\begin{equation}\label{ex:eq:relative-strict-minimum}
 0<\Flat(T-T_i)<\rho_i
 \quad\Longrightarrow\quad \Mass(T)>4\pi.
\end{equation}
Fix $i$ and let $E_i=[0,t_i]\times S^2$, so
$\dr\llbracket E_i\rrbracket=T_i$. Write
$P(E;B)=\HH^2(\partial^*E\cap B\cap\Int M)$ and $P(E)=P(E;M)$.
For small $\Flat(T-T_i)$, Lemma~\ref{inp:compactness} gives a filling
$Q=\llbracket A\rrbracket$ with $\Mass(Q)\le C_M\Flat(T-T_i)$,
represented by a finite-perimeter set
\cite[Sections 4.2 and 4.5]{Federer1969}.
Then $E=E_i\triangle A$ represents $T$, with
\begin{equation}\label{ex:eq:relative-set-comparison}
 \Mass(T)=P(E),\quad
 v:=\Vol(E\triangle E_i)\le C_M\Flat(T-T_i).
\end{equation}
Indeed, relative boundaries add under symmetric difference, so
$\dr\llbracket E\rrbracket=T_i+(T-T_i)=T$. The canonical mass of
this boundary is $P(E)$. If $T\ne T_i$, then $v>0$; it therefore
suffices to prove $P(E)>4\pi$ whenever $v$ is small and positive.

Put $z=t-t_i$ and write $U_b=\{|z|<b\}$. By \eqref{ex:eq:slice-derivatives}, we can choose $a,c>0$ such that $\overline{U_{2a}}\subset\Int M$ and
$f(t)^2\ge1+2cz^2$ on $U_{2a}$.
Choose a smooth function $0\le\eta\le1$, equal to one on $U_a$ and supported in $U_{2a}$, and define
\[
 X=\eta\frac{1+cz^2}{f(t)^2}\partial_t
\]
on $U_{2a}$, extending it by zero. We have $|X|\le1$, $X=\partial_t$ on $\Sigma_i$, and, for some $\gamma>0$, $|X|\le1-\gamma$ on $M\setminus U_{a/2}$. On $U_a$, direct differentiation gives
\begin{equation}\label{ex:eq:radial-divergence}
 \operatorname{div}_gX
 =\frac1{f(t)^2}\frac{d}{dt}(1+cz^2)
 =\frac{2cz}{f(t)^2}.
\end{equation}
All constants in the following estimates may depend on $i$.

Let $\nu_E$ be the measure-theoretic outward unit normal of $E$, and set
\[
 L=P(E;M\setminus U_{a/2}),\quad
 W=\int_{A\cap U_a}|z|\,d\Vol.
\]
Since $X$ is supported in $\Int M$ and
$\int_{E_i}\diver_gX=4\pi$, Gauss--Green gives
\begin{equation}\label{ex:eq:perimeter-comparison}
 \begin{aligned}
 P(E)-4\pi
 &=\int_{\partial^*E\cap\Int M}(1-\langle X,\nu_E\rangle)\,d\HH^2
   +\int_M(\mathbf1_E-\mathbf1_{E_i})\diver_gX\,d\Vol\\
 &\ge\gamma L+c_0W-C\Vol(A\setminus U_a).
 \end{aligned}
\end{equation}
Here $c_0,C>0$:
on $A\cap U_a$, $\mathbf1_E-\mathbf1_{E_i}$ has the sign of $z$,
so \eqref{ex:eq:radial-divergence} gives the $W$ term.

To estimate the last term, coarea gives
\[
 \int_{a/2}^{a}\HH^2(A\cap\partial U_b)\,db
 =\Vol\bigl(A\cap(U_a\setminus U_{a/2})\bigr)
 \le\frac{2W}{a}.
\]
Choose a regular $b\in(a/2,a)$ with
$\HH^2(A\cap\partial U_b)\le4W/a^2$, and set $F=A\setminus U_b$.
Almost everywhere, $A=M\setminus E$ below $t_i$ and $A=E$ above it;
their perimeters thus agree outside $U_{a/2}$. Restriction yields
$P(F)\le L+4W/a^2$.
The Euclidean ball inequality \cite[Section 5.6.2, Theorem 2]{EG1992}
and uniform metric comparison give
$\Vol(F)^{2/3}\le C P(F)$ for $\Vol(F)\le\Vol(M)/2$.
Since $F\subset A$ and $\Vol(A)=v$, we may multiply this inequality
by $\Vol(F)^{1/3}\le v^{1/3}$. Hence, for $v\le\Vol(M)/2$,
\begin{equation}\label{ex:eq:exterior-volume}
 \Vol(A\setminus U_a)\le\Vol(F)
 \le C v^{1/3}P(F)\le C v^{1/3}(L+W).
\end{equation}
Combining \eqref{ex:eq:perimeter-comparison} and \eqref{ex:eq:exterior-volume}, and decreasing the upper bound for $v$, gives
\[
 P(E)-4\pi\ge\frac{\gamma}{2}L+\frac{c_0}{2}W.
\]
If $L=W=0$, \eqref{ex:eq:exterior-volume} gives
$\Vol(A\setminus U_a)=0$, while $W=0$ gives $\Vol(A\cap U_a)=0$.
For the latter implication, $|z|>0$ away from the slice $\{z=0\}$,
which has zero three-dimensional volume.
Thus $v=0$, proving strict inequality for small $v>0$ and, by
\eqref{ex:eq:relative-set-comparison}, \eqref{ex:eq:relative-strict-minimum}.

\noindent\textbf{Almost minimization.}

Fix $i$ and choose $0<r_i<\rho_i/3$. By relative flat compactness and lower semicontinuity of mass, the set
\[
 \mathcal K_i=\{T:r_i\le\Flat(T-T_i)\le2r_i,
                         \ \Mass(T)\le4\pi+1\}
\]
is compact by Lemma~\ref{inp:compactness}. If $\mathcal K_i$ is nonempty, \eqref{ex:eq:relative-strict-minimum} gives
$\beta_i=\min_{T\in\mathcal K_i}\Mass(T)-4\pi>0$;
the minimum is attained by lower semicontinuity, and its minimizer
lies at positive flat distance less than $\rho_i$ from $T_i$.
If $\mathcal K_i=\varnothing$, set $\beta_i=1$.
Choose $0<\delta_i<\frac12\min\{r_i,\beta_i,1\}$. Consider any sequence $T_i=T^0,\ldots,T^N$ with relative flat steps at most $\delta_i$ and masses at most $4\pi+\delta_i$. If $j$ is the first index with $\Flat(T^j-T_i)\ge r_i$, then
\[
 r_i\le\Flat(T^j-T_i)<r_i+\delta_i<2r_i,
 \quad \Mass(T^j)\le4\pi+\delta_i<4\pi+\beta_i.
\]
This contradicts the definition of $\beta_i$ (or the emptiness of
$\mathcal K_i$). Hence every term stays within distance $r_i$, and
\eqref{ex:eq:relative-strict-minimum} yields
$\Mass(T^N)\ge4\pi$. Consequently
$T_i\in\AMset(M;\eps,\delta_i;\Flat)$ for every $\eps>0$.

Since $|T_i|\rightharpoonup V=|\partial M|$, choose
$\dV(|T_i|,V)<\eps$ and use its $\delta_i$ in
Definition~\ref{def:almostminimizing}. This proves almost minimization
in $M$ and every relative annulus, without a uniform lower bound on
$\delta_i$. The divergence theorem on $\partial M$ gives free-boundary
stationarity; its area gives \eqref{ex:eq:boundary-concentration}.
\end{proof}

The relative currents and their varifolds have different limits in this example:
\[
 \Flat(T_i)\le\Vol([0,t_i]\times S^2)
 =4\pi\int_0^{t_i}f(t)^2\,dt\longrightarrow0,
 \quad |T_i|\rightharpoonup|\partial M|\ne0.
\]
Proposition~\ref{prop:boundary-concentration} thus uses the relative
almost minimizing condition of Theorem~\ref{thm:kflat};
Example~\ref{ex:spherical-boundary} compares absolute cycles.

\clearpage
\renewcommand{\bibfont}{\small}
\setlength{\bibsep}{2pt}

\end{document}